\documentclass[10pt]{article}
\usepackage{latexsym}
\usepackage{amssymb}
\usepackage{amsmath,amscd}
\usepackage{amsthm}
\usepackage{mathrsfs}
\usepackage{color}
\usepackage{graphicx}
\allowdisplaybreaks

\newtheorem{theorem}{Theorem}[section]
\newtheorem{lemma}{Lemma}[section]
\newtheorem{proposition}{Proposition}[section]
\newtheorem{definition}{Definition}[section]
\newtheorem{remark}{Remark}[section]

\newtheorem{example}{Example}[section]

\begin{document}

\title{\bf Conley index for random dynamical systems\footnote{This
work is partially supported by `985' project of Jilin
University and graduate innovation lab of Jilin University.}}

\author{Zhenxin Liu\\
{\small College of Mathematics, Jilin University,  Changchun 130012,
P.R. China}\\
{\small zxliu@email.jlu.edu.cn}}

\date{}
\maketitle

\begin{abstract}
Conley index theory is a very powerful tool in the study of
dynamical systems, differential equations and bifurcation theory. In
this paper, we make an attempt to generalize the Conley index to
discrete random dynamical systems. And we mainly follow the Conley
index for maps given by Franks and Richeson in \cite{Fra}.
Furthermore, we simply discuss the relations of isolated invariant
sets between time-continuous random dynamical systems and the
corresponding time-$h$ maps. For applications we give several
examples to illustrate our results.
\end{abstract}

\section{Introduction}

One important aspect of the qualitative analysis of differential
equations and dynamical systems is the study of asymptotic,
long-term behavior of solutions/orbits. Hence much of dynamical
systems involves the study of the existence and structure of
invariant sets. The Conley index theory (see \cite{Con,CZ,Sal}),
developed by Conley and his students, has been a very powerful tool
in the study of dynamical systems, differential equations and
bifurcation theory. The Conley index theory was developed for flows
on compact spaces at first, then Rybakowski \cite{Ry} extended the
theory to semiflows on noncompact spaces. After that, a natural
question is to find an appropriate generalization to discrete
dynamical systems (maps). Robbin and Salamon \cite{Rob} presented
the first version of the Conley index for maps, which was defined in
terms of shape theory. Later Mrozek \cite{Mro} offered a
cohomological definition of Conley index for maps based on the Leray
reduction. Szymczak \cite{Szy} constructed the homotopy Conley index
as a functor into an abstract category and showed that any other
Conley index can be factorized through his homotopy Conley index.
Recently, a more intuitive definition of the discrete Conley index
was given by Franks and Richeson \cite{Fra}, which gives an
accessible and intuitive development of a Conley index for isolated
invariant sets of any continuous map defined on a locally compact
metric space. They also showed that their definition of the discrete
Conley index based on shift equivalences is identical to Szymczak's
categorical definition of the Conley index.

Random dynamical systems arise in the modeling of many phenomena in
physics, biology, economics, climatology, etc and the random effects
often reflect intrinsic properties of these phenomena rather than
just to compensate for the defects in deterministic models. The
history of study of random dynamical systems goes back to Ulam and
von Neumann \cite{Ula} and it has flourished since the 1980s due to
the discovery that the solutions of stochastic ordinary differential
equations yield a cocycle over a metric dynamical system which
models randomness, i.e. a random dynamical system. In this paper we
make an attempt to obtain the Conley index for random dynamical
systems. The definition of Conley index given by Franks and Richeson
\cite{Fra} is relative simple and unlike index pairs, the filtration
pairs they defined are robust under small $C^0$ perturbations of the
map $f$, which is a wonderful property to obtain the continuation
property of Conley index immediately. Hence we will follow
\cite{Fra} to obtain the Conley index for discrete random dynamical
systems. Such Conley index is very useful in studying random
homeomorphisms, discretized random differential equations etc, see
Section 8 for several simple applications. Furthermore we also
discuss the isolated invariant set for time-continuous random
dynamical systems and its relation with the isolated invariant set
for the time-$h$ maps. This indicates that we can obtain some
information of time-continuous random dynamical system by studying
its time-$h$ maps.

The paper is organized as follows. In Section 2, we recall some
basic definitions and results for random dynamical systems. In
Section 3, we give the definitions of random isolating neighborhood
and random isolated invariant set and compare the random isolated
invariant set and random omega-limit set. In Section 4, we prove
that for any random neighborhood of a random isolated invariant set,
there is a random filtration pair for it in this neighborhood. In
Section 5, we introduce the random shift equivalence and show that
for any two random filtration pairs of a random isolated invariant
set, the corresponding two random pointed spaces with their random
pointed space maps are random shift equivalent. In Section 6 we give
the definition of random Conley index for random isolated invariant
sets. In Section 7 we simply discuss the relation of isolated
invariant sets between time-continuous random dynamical system and
the discrete one generated by its time-$h$ maps. And at last for
applications we give several simple examples to illustrate our
results in Section 8.

\section{Random dynamical systems}

In this section, we will give some preliminary definitions and
propositions for the later use. Firstly we give the definition of
continuous random dynamical systems (cf Arnold \cite{Ar1}).

\begin{definition}\rm
Let $X$ be a metric space with a metric $d_X$. A continuous {\em
random dynamical system (RDS)}, shortly denoted by
$\phi$, consists of two ingredients: \\
(i) A model of the noise, namely a metric dynamical system $(\Omega,
\mathscr F, \mathbb P, (\theta_t)_{t\in \mathbb T})$, where
$(\Omega, \mathscr F, \mathbb P)$ is a probability space and $(t,
\omega)\mapsto \theta_t\omega$ is a measurable flow which leaves
$\mathbb P$ invariant, i.e. $\theta_t\mathbb P=\mathbb P$ for all
$t\in \mathbb T$.\\
(ii) A model of the system perturbed by noise, namely a cocycle
$\phi$ over $\theta$, i.e. a measurable mapping $\phi: \mathbb
T\times \Omega\times X \rightarrow X,
(t,\omega,x)\mapsto\phi(t,\omega,x)$, such that
$(t,x)\mapsto\phi(t,\omega,x)$ is continuous for all
$\omega\in\Omega$ and the family
$\phi(t,\omega,\cdot)=\phi(t,\omega):X\rightarrow X$ of random
self-mappings of $X$ satisfies the cocycle property:
\begin{equation}\label{phi}
\phi(0,\omega)={\rm id}_X, \phi(t+s,\omega)=\phi(t,\theta_s
\omega)\circ\phi(s,\omega)\quad {\rm for ~all}\quad t,s\in\mathbb
T,\omega\in\Omega.
\end{equation}
In this definition, $\mathbb T=\mathbb Z$ or $\mathbb R$.
\end{definition}

In this paper, we mainly consider discrete random dynamical systems,
i.e. $\mathbb T=\mathbb Z$. We will use $\theta:=\theta_1$ to denote
the time one map of $\theta_n$. It follows from (\ref{phi}) that
$\phi(n,\omega)$ is a homeomorphism of $X$, and the fact
\[\phi^{-1}(n,\omega)=\phi(-n,\theta_n\omega)\]
is very useful in the following.

Assume $\phi$ is a discrete random dynamical system and $\varphi$ is
the time-one map of $\phi$, i.e.
$\varphi(\omega)=\phi(1,\omega):X\rightarrow X$, then we call
$\varphi$ the {\em random homeomorphism} determined by $\phi$. On
the other hand, assume $\varphi$ is a random homeomorphism, then it
generates a discrete RDS $\phi(n,\omega,x)$ in the following way:
\[
\phi(n,\omega)=\left\{
                  \begin{array}{ll}
                  \varphi(\theta_n\omega)\circ\cdots\circ\varphi(\omega) , & n>0, \\
                   {\rm Id}_X , & n=0,\\
                   \varphi^{-1}(\theta_n\omega)\circ\cdots\circ\varphi^{-1}(\theta_{-1}\omega), & n<0.
                  \end{array}
                \right.
\]
Hence we will identify a random homeomorphism with the discrete RDS
generated by it, which will not cause confusions. For a random
homeomorphism $\varphi$, we will use $\varphi^k$ to denote the
k-times iteration of $\varphi$, i.e.
$\varphi^k(\omega,\cdot):=\phi(k,\omega,\cdot)$, where $\phi$ is the
discrete RDS generated by $\varphi$.

Below any mapping from $\Omega$ into the collection of all subsets
of $X$ is said to be a {\it multifunction} (or a set valued mapping)
from $\Omega$ into X. We now give the definition of random set,
which is a fundamental concept for RDS.

\begin{definition}\rm
Let $X$ be a metric space with a metric $d_X$. The multifunction
$\omega\mapsto D(\omega)$ taking values in the closed/compact
subsets of $X$ is said to be a {\em random closed/compact set} if
the mapping $\omega\mapsto {\rm dist}_X(x,D(\omega))$ is measurable
for any $x \in X$, where ${\rm dist}_X(x,B):=\inf_{y\in B}d_X(x,y)$.
The multifunction $\omega\mapsto U(\omega)$ taking values in the
open subsets of $X$ is said to be a {\em random open set} if
$\omega\mapsto U^c(\omega)$ is a random closed set, where $U^c$
denotes the complement of $U$.
\end{definition}

\begin{definition}\rm
A random set $D(\omega)$ is said to be {\em forward invariant} under
the RDS $\phi$ if $\phi(t,\omega)D(\omega)\subset D(\theta_t\omega)$
for all $t\in\mathbb T_+$; It is said to be {\em backward invariant}
if $\phi(t,\omega)D(\omega)\supset D(\theta_t\omega)$ for all
$t\in\mathbb T_+$; It is said to be {\em invariant} if
$\phi(t,\omega)D(\omega)=D(\theta_t\omega)$ for all $t\in\mathbb T$.
\end{definition}

Throughout the paper, we will assume that $(X,d_X)$ is a locally
compact {\em Polish space}, i.e. a separable complete metric space.
For a random variable $T(\omega)$, we call $T(\omega)>0$ if it holds
almost surely. We also call a multifunction $D(\omega)$ measurable
for convenience if the mapping $\omega\mapsto {\rm
dist}_X(x,D(\omega))$ is measurable for any $x \in X$.

Now we enumerate some basic results about random sets in the
following propositions, for details the reader can refer to Castaing
and Valadier \cite{Cas} and Arnold \cite{Ar1} for instance.
\begin{proposition}\label{set}
Let X be a Polish space, then the following assertions hold: \\
(i) if D is a random closed set, then so is the  closure of $D^c$;\\
(ii) if D is a random open set, then the closure ${\rm cl}D$ of $D$
is a random closed set;\\
(iii) if D is a random closed set, then {\rm int}$D$, the interior
of $D$, is a random open set; \\
(iv) if $\{D_n, n \in\mathbb N\}$ is a sequence of random closed
sets and there exists $n_0 \in\mathbb N$ such that $D_{n_0}$ is a
random compact set, then
$\bigcap_{n\in\mathbb N}D_n$ is a random compact set;\\
(v) if $f:\Omega\times X\rightarrow X$ is a function such that
$f(\omega,\cdot)$ is continuous for all $\omega$ and $f(\cdot,x)$ is
measurable for all $x$, then $\omega\mapsto f(\omega,D(\omega))$ is
a random compact set provided that $D(\omega)$ is a random compact
set.
\end{proposition}

The following measurable selection theorem is frequently used for
our purpose, although we do not always mention it when we do use it.
For the proof, we refer to \cite{Ar1,Cas} for details.
\begin{proposition} {\rm(Measurable Selection
Theorem)}. Let a multifunction $\omega\mapsto D(\omega)$ take values
in the subspace of closed non-void subsets of a Polish space $X$.
Then $D(\omega)$ is a random closed set if and only if there exists
a sequence $\{v_n : n \in\mathbb N\}$ of measurable maps $v_n
:\Omega\rightarrow X$ such that
\[
 v_n(\omega)\in D(\omega)\quad and
\quad D(\omega) ={\rm cl}(\bigcup{\{v_n(\omega),n\in\mathbb
N\}})\quad for~all \quad \omega\in\Omega.
\]
In particular if $D(\omega)$ is a random closed set, then there
exists a measurable selection, i.e. a measurable map $v
:\Omega\rightarrow X$  such that $v(\omega)\in D(\omega)$ for all
$\omega\in\Omega$.
\end{proposition}

The following proposition comes from \cite{Rud} and will be used
later, which gives a relation between an $\bar{\mathscr
F}^\nu$-measurable function and an $\mathscr F$-measurable one.
\begin{proposition}\label{rud}
Assume $\nu$ is a positive measure on the measurable space
$(X,\mathscr F)$. Denote $\bar{\mathscr F}^\nu$ the completion of
the $\sigma$-algebra $\mathscr F$ with respect to the measure $\nu$.
If $f$ is an $\bar{\mathscr F}^\nu$-measurable function, then there
exists an $\mathscr F$-measurable function $g$ such that $f=g$ a.e.
$[\nu]$.
\end{proposition}

\section{Random isolated invariant sets and omega-limit sets}

\begin{definition}\rm
For given two random sets $D(\omega),A(\omega)$, we say $A(\omega)$
{\em attracts} $D(\omega)$ if
\[
\lim_{n\rightarrow\infty}d(\phi(n,\theta_{-n}\omega)D(\theta_{-n}\omega)|A(\omega))=0
\]
holds almost surely, where $d(A|B)$ stands for the Hausdorff
semi-metric between two sets $A,B$, i.e. $d(A|B):={\rm sup}_{x\in
A}{\rm inf}_{y\in B}d(x,y)$.
\end{definition}

For any given random set $D(\omega)$, we denote $\Omega_D(\omega)$
the {\em omega-limit set of $D(\omega)$}, which is defined as
follows:
\[
\Omega_D(\omega):=\bigcap_{n\ge 0}\overline{\bigcup_{k\ge
n}\phi(k,\theta_{-k}\omega)D(\theta_{-k}\omega)}.
\]
It is well-known that if a nonvoid  random set $D(\omega)$ is
attracted by a random compact set $A(\omega)$, then
$\Omega_D(\omega)\neq\emptyset$ almost surely and it is invariant.
Also it is known that for two random sets $D_1(\omega)$,
$D_2(\omega)$, if they are attracted by a random compact set
$A(\omega)$, then we have
\[
\Omega_{D_1\cup
D_2}(\omega)=\Omega_{D_1}(\omega)\cup\Omega_{D_2}(\omega).
\]
This important fact will be used later, e.g. in Lemma \ref{lem4}.
(In fact, this conclusion also holds if $D_1(\omega)$, $D_2(\omega)$
are not attracted by a random compact set.)

\begin{definition}\label{def}\rm
A random compact set $N(\omega)$ is called a {\em random isolating
neighborhood} if it satisfies
\[
{\rm Inv}N(\omega)\subset {\rm int}N(\omega),
\]
where ${\rm int}N(\omega)$ denotes the interior of $N(\omega)$ and
\[
{\rm Inv}N(\omega)=\{x\in N(\omega)|~\phi(n,\omega,x)\in
N(\theta_n\omega),\forall n\in\mathbb Z\}.
\]
Correspondingly we call $S(\omega)$ a {\em random isolated invariant
set} if there exists a random isolating neighborhood $N(\omega)$
such that $S(\omega)={\rm Inv}N(\omega)$. A random compact set
$N(\omega)$ is called a {\em random isolating block} if it satisfies
\[
\varphi(\theta_{-1}\omega,N(\theta_{-1}\omega))\cap
N(\omega)\cap\varphi^{-1}(\theta\omega,N(\theta\omega))\subset{\rm
int}N(\omega).
\]
\end{definition}

\begin{definition}\rm
For a random set $N(\omega)$ we define the {\em exit set of
$N(\omega)$} to be
\[
N^-(\omega):=\{x\in N(\omega)|~\varphi(\omega,x)\not\in {\rm
int}N(\theta\omega)\}.
\]
\end{definition}

\begin{remark}\label{r1}\rm
It is obvious that a random isolating block is a random isolating
neighborhood, but the converse is not true. It is easy to see that
for a random isolating neighborhood $N(\omega)$, its corresponding
random isolated invariant set $S(\omega)$ can be characterized by
\[
S(\omega)={\rm Inv}N(\omega)=\bigcap_{n\in\mathbb
Z}\phi(-n,\theta_n\omega)N(\theta_n\omega)
\]
and the exit set $N^-(\omega)$  of a random set $N(\omega)$ can be
characterized  by
\begin{equation}\label{N-}
N^-(\omega)=N(\omega)\cap[\varphi^{-1}(\theta\omega,{\rm
int}^cN(\theta\omega))],
\end{equation}
where ${\rm int}^cN(\omega)$ denotes the complement of the interior
of $N(\omega)$. It is also obvious that
\[
\varphi(\omega,N^-(\omega))\cap {\rm
int}N(\theta\omega)=\emptyset,~\varphi(\omega,N(\omega)\backslash
N^-(\omega))\subset {\rm int}N(\theta\omega).
\]
These simple facts will be used in the sequel.
\end{remark}

The following remark is useful in our proof, see e.g. Lemmas
\ref{lem3} and \ref{lem4}.
\begin{remark}\label{r2}\rm
(i) Assume $N(\omega)$ is a random compact set, it is obvious that
${\rm Inv}N(\omega)\subset\Omega_N(\omega)$ and it is easy to see
that ${\rm Inv}N(\omega)=\Omega_N(\omega)$ if and only if
$\Omega_N(\omega)\subset N(\omega)$. In particular, when $N(\omega)$
is forward invariant, we have $\Omega_N(\omega)\subset N(\omega)$.\\
(ii) Unlike omega-limit sets, in general we only have
\[
{\rm Inv}D_1(\omega)\cup{\rm Inv}D_2(\omega)\subset{\rm Inv}(D_1\cup
D_2)(\omega),
\]
but not ${\rm Inv}D_1(\omega)\cup{\rm Inv}D_2(\omega)={\rm
Inv}(D_1\cup D_2)(\omega)$.
\end{remark}

\section{Random filtration pair}

\begin{definition}\rm
Assume $N(\omega)$ is a random isolating neighborhood,
$L(\omega)\subset N(\omega)$ is a random compact set and $S(\omega)$
is the random isolated invariant set inside $N(\omega)$. We also
assume $N(\omega)={\rm cl(int}N(\omega))$, $L(\omega)={\rm
cl(int}L(\omega))$. We call $(N(\omega),L(\omega))$ is a {\em random
filtration pair for $S(\omega)$} if the following holds:
\begin{itemize}
    \item ${\rm cl}(N(\omega)\backslash L(\omega))$ is a random
    isolating neighborhood of $S(\omega)$;
    \item $L(\omega)$ is a random neighborhood of $N^-(\omega)$ in
    $N(\omega)$;
    \item $\varphi(\omega,L(\omega))\cap{\rm cl}(N(\theta\omega)\backslash
    L(\theta\omega))=\emptyset$.
\end{itemize}
\end{definition}

Similar to \cite{Liu}, for a given random variable
$\epsilon(\omega)>0$ we define a {\em random
$\epsilon(\omega)$-chain of length $n$} to be $n+1$ random variables
$x_0(\omega)$, $x_1(\omega)$, $\cdots$, $x_n(\omega)$ satisfying
that
\[
d_X(\varphi(\theta_{-1}\omega,x_i(\theta_{-1}\omega)),x_{i+1}(\omega))<\epsilon(\omega),
\]
for $i=0,1,\cdots,n-1$.

\begin{definition}\rm
Assume $N(\omega)$ is a random isolating neighborhood and
$S(\omega)$ is the random isolated invariant set inside $N(\omega)$.
For any random variable $\epsilon(\omega)>0$, we define the {\em
random $\epsilon(\omega)$-chain neighborhood of $S(\omega)$}
relative to $N(\omega)$, $C_\epsilon(N,S)(\omega)$, to be the union
of the random variables $x(\omega)$ such that
$\{x_i(\omega)\}_{-l}^{m}\subset N(\omega)$ is a random
$\epsilon(\omega)$-chain satisfying that
$x_{-l}(\omega),x_m(\omega)\in S(\omega)$ and
$x_0(\omega)=x(\omega)$.
\end{definition}

For the properties of the random $\epsilon(\omega)$-chain
neighborhood of $S(\omega)$ relative to $N(\omega)$, see the
following several lemmas.

\begin{lemma}
Assume $N(\omega)$ is a random isolating neighborhood and
$S(\omega)$ is the random isolated invariant set inside $N(\omega)$.
Then $C_\epsilon(N,S)(\omega)$ is a random open set relative to
$N(\omega)$.
\end{lemma}
{\it Proof.} Define
\begin{align*}
& U_0(\omega):=S(\omega),\\
&
U_1(\omega):=B_{\epsilon(\omega)}(\varphi(\theta_{-1}\omega,S(\theta_{-1}\omega)))
(=B_{\epsilon(\omega)}(S(\omega)))\\
&
U_2(\omega):=B_{\epsilon(\omega)}(\varphi(\theta_{-1}\omega,U_1(\theta_{-1}\omega)))\\
& \cdots\cdots\\
&
U_n(\omega):=B_{\epsilon(\omega)}(\varphi(\theta_{-1}\omega,U_{n-1}(\theta_{-1}\omega)))\\
& \cdots\cdots
\end{align*}
and define
\begin{align*}
& V_0(\omega):=S(\omega),\\
&
V_1(\omega):=\varphi^{-1}(\theta\omega,B_{\epsilon(\theta\omega)}(V_0(\theta\omega)))\\
&
\cdots\cdots\\
&
V_n(\omega):=\varphi^{-1}(\theta\omega,B_{\epsilon(\theta\omega)}(V_{n-1}(\theta\omega)))\\
& \cdots\cdots,
\end{align*}
where $B_r(A)$ stands for the open $r$-neighborhood of the set $A$.
Denote
\[
U(\omega)=\bigcup_{i=0}^{\infty}U_i(\omega),~V(\omega)=\bigcup_{i=0}^{\infty}V_i(\omega),
\]
then by the definition of $C_\epsilon(N,S)(\omega)$ we have
\begin{equation*}
C_\epsilon(N,S)(\omega)=U(\omega)\cap V(\omega)\cap N(\omega).
\end{equation*}
By the fact $\varphi(\omega,\cdot)$ is homeomorphism for
$\forall\omega$, we obtain that each $U_i(\omega),V_i(\omega)$,
$i=1,2,\ldots$ is a random open set. Then we obtain that
$U(\omega)$, $V(\omega)$ are random open sets. It is easy to see
that $U^c(\omega)\cup V^c(\omega)$ is measurable, i.e. the
complement of $U(\omega)\cap V(\omega)$ is a random closed set. So
$U(\omega)\cap V(\omega)$ is a random open set. Hence
$C_\epsilon(N,S)(\omega)$ is a random open set relative to
$N(\omega)$. \hfill $\Box$

\begin{lemma}\label{S1}
Assume $N(\omega)$ is a random isolating neighborhood and
$S(\omega)$ is the random isolated invariant set inside $N(\omega)$.
Then $S(\omega)$ can be characterized by
\begin{equation}\label{S}
S(\omega)=\bigcap\{{\rm
cl}(C_\epsilon(N,S)(\omega))|~\epsilon(\omega)>0\}.
\end{equation}
\end{lemma}
{\it Proof.} The idea of the proof is originated from Easton
\cite{Eas}. It is obvious that the left hand of (\ref{S}) is a
subset of right hand, so we only need to verify that the converse
inclusion is true almost surely. For arbitrary random variable
$x(\omega)\in\bigcap\{{\rm
cl}(C_\epsilon(N,S)(\omega))|~\epsilon(\omega)>0\}$, and
$\forall\epsilon(\omega)>0$, there exists a random
$\epsilon(\omega)$-chain in $N(\omega)$ from $x(\omega)$ to
$S(\omega)$ by the definition of $C_\epsilon(N,S)(\omega)$. Let
$\epsilon(\omega)\rightarrow 0$, by the invariance of $S(\omega)$
and the measure preserving of $\theta_n$, we obtain that the entire
forward orbit of $x(\omega)$ is in $N(\omega)$ almost surely, i.e.
$\varphi^n(\omega,x(\omega))\in N(\theta_n\omega)$ for $\forall n\in
\mathbb N$. Similarly the entire backward orbit of $x(\omega)$ is
also in $N(\omega)$ almost surely. Hence we have obtained
\[
\varphi^n(\omega,x(\omega))\in N(\theta_n\omega),\forall n\in
\mathbb Z,
\]
i.e. $x(\omega)\in S(\omega)$ almost surely by the definition of
$S(\omega)$. Thus we have verified that the the right hand of
(\ref{S}) is a subset of $S(\omega)$ almost surely. \hfill $\Box$

\begin{lemma}\label{lem1}
Assume $S(\omega)$ is a random isolated invariant set with a random
isolating neighborhood $N(\omega)$ and $W(\omega)$ is an arbitrary
random neighborhood of $S(\omega)$ in $N(\omega)$. Then there exists
a random variable $\epsilon(\omega)>0$ such that
$C_\epsilon(N,S)(\omega)\subset W(\omega)$.
\end{lemma}
{\it Proof.} By the invariance of $S(\omega)$, it is easy to see
that for a given $k\in\mathbb N$, there exists a sufficiently small
random variable $\tilde\epsilon(\omega)$ such that all the random
$\tilde\epsilon(\omega)$-chains with length not greater than $k$
from $S(\omega)$ to $S(\omega)$ must be in $W(\omega)$. Denote
$\epsilon_n(\omega)=\frac1{2^n}\tilde\epsilon(\omega)$, then by a
similar argument to that of Lemma \ref{S1} we obtain that
\[
S(\omega)=\bigcap_{n\in\mathbb N}\{{\rm
cl}(C_{\epsilon_n}(N,S)(\omega))\}.
\]
For simplicity, we denote $S_n(\omega)={\rm
cl}(C_{\epsilon_n}(N,S)(\omega))$. It is obvious that
$S_n(\omega)\supset S_{n+1}(\omega)$, hence we have
\[
S(\omega)=\bigcap_{n\in\mathbb
N}S_n(\omega)=\lim_{n\rightarrow\infty}S_n(\omega).
\]
Therefore, $\exists\tilde n=\tilde n(\omega)$ such that
$S_n(\omega)\subset W(\omega)$ whenever $n\ge\tilde n$. In fact, if
we take
\[
\tilde n(\omega):=\inf\{n\in\mathbb N|~d(S_n(\omega)|W(\omega))=0\},
\]
then, similar to Lemma 3.5 of \cite{Liu} and by Proposition
\ref{rud}, we obtain that $\tilde n(\omega)$ is measurable. Define
$\epsilon(\omega):=\epsilon_{\tilde n}(\omega)$, then
$\epsilon(\omega)$ is measurable and by its definition we have
$C_{\epsilon}(N,S)(\omega)\subset W(\omega)$. This completes the
proof of the lemma. \hfill $\Box$

\begin{lemma}\label{lem2}
Assume $S(\omega)$ is a random isolated invariant set with a random
isolating neighborhood $N(\omega)$ and $\epsilon(\omega)>0$ is
sufficiently small, then ${\rm cl}(C_{\epsilon}(N,S)(\omega))$ is a
random isolating block.
\end{lemma}
{\it Proof.} For $\forall x(\omega)\in C_{\epsilon}(N,S)(\omega)$,
by the definition of $C_{\epsilon}(N,S)(\omega)$, we have the
following holds:
\[
B_{\epsilon(\omega)}(\varphi(\theta_{-1}\omega,x(\theta_{-1}\omega)))\cap
N(\omega)\cap
\varphi^{-1}(\theta\omega,B_{\epsilon(\theta\omega)}(x(\theta\omega)))\subset
C_{\epsilon}(N,S)(\omega).
\]
By the fact that $\varphi(\omega,\cdot)$ is a homeomorphism on $X$,
there exists a $\delta(\omega)>0$ such that
\[
B_{\delta(\omega)}(\varphi(\theta_{-1}\omega,x(\theta_{-1}\omega)))\cap
N(\omega)\cap
B_{\delta(\omega)}(\varphi^{-1}(\theta\omega,x(\theta\omega)))
\subset C_{\epsilon}(N,S)(\omega).
\]
Since $x(\omega)\in C_{\epsilon}(N,S)(\omega)$ is arbitrary, we
obtain that
\begin{align*}
&\varphi(\theta_{-1}\omega,{\rm
cl}(C_{\epsilon}(N,S)(\theta_{-1}\omega)))\cap N(\omega)\cap
\varphi^{-1}(\theta\omega,{\rm
cl}(C_{\epsilon}(N,S)(\theta\omega)))\\
&\subset C_{\epsilon}(N,S)(\omega)\subset{\rm int}({\rm
cl}(C_{\epsilon}(N,S)(\omega))),
\end{align*}
i.e. ${\rm cl}(C_{\epsilon}(N,S)(\omega))$ is a random isolating
block. This completes the proof of the lemma. \hfill $\Box$

Assume $\varphi(\omega,\cdot),\psi(\omega,\cdot)$ are two random
homeomorphisms and $N(\omega)$ is a random compact set. We define
\begin{align*}
d_{\varphi,\psi}^N(\omega)&:=\sup_{x\in
N(\theta_{-1}\omega)}d_X(\varphi(\theta_{-1}\omega,x),\psi(\theta_{-1}\omega,x))\\
&\qquad+\sup_{y\in
N(\theta\omega)}d_X(\varphi^{-1}(\theta\omega,y),\psi^{-1}(\theta\omega,y)).
\end{align*}
Then by the measurable selection theorem we have
\begin{align}\label{norm}
d_{\varphi,\psi}^N(\omega)&=\sup_{i=1}^\infty
d_X(\varphi(\theta_{-1}\omega,x_i(\theta_{-1}\omega)),
\psi(\theta_{-1}\omega,x_i(\theta_{-1}\omega)))\nonumber\\
&\qquad+\sup_{j=1}^\infty
d_X(\varphi^{-1}(\theta\omega,x_j(\theta\omega)),\psi^{-1}(\theta\omega,x_j(\theta\omega))),
\end{align}
where $\{x_n(\omega)\}_{n=1}^{\infty}$ is a family of countable
dense random variables of $N(\omega)$. Hence
$d_{\varphi,\psi}^N(\omega)$ is measurable. It is easy to see that
for any given $\omega\in\Omega$, $d_{\varphi,\psi}^N(\omega)$ is a
metric defined on the space of random homeomorphisms whose
definition domain contains $N(\omega)$ and we call it ``random
metric". In the following, we say a family of random homeomorphisms
$\psi_n$ approximate $\psi$ in the random $C^0$ topology if
$d_{\psi_n,\psi}^N(\omega)$ converge to $0$ when
$n\rightarrow\infty$.

The following theorem states that the random filtration pairs for a
random isolated invariant set are robust under small random $C^0$
perturbations of the discrete random dynamical system $\varphi$,
through which we can obtain the continuation property of random
Conley index easily.

\begin{theorem}\label{th1}
Assume $N(\omega)$ is a random isolating block and $L(\omega)$ is a
sufficiently small random compact neighborhood of $N^-(\omega)$ in
$N(\omega)$, then $(N(\omega),L(\omega))$ is a random filtration
pair for ${\rm Inv}N(\omega)$. Moreover, there exists a random
$C^0$-neighborhood of $\varphi$ such that for any random
homeomorphism $\psi$ in this neighborhood $S_\psi(\omega):={\rm
Inv}(N(\omega)\backslash L(\omega),\psi)$ is a random isolated
invariant set for $\psi$ and $(N(\omega),L(\omega))$ is a random
filtration pair for $S_\psi(\omega)$.
\end{theorem}
{\it Proof.} (1) By the definition of $N^-(\omega)$, we have
\begin{equation}\label{1}
\varphi(\omega,N^-(\omega))\subset {\rm int}^cN(\theta\omega).
\end{equation}
Denote $S(\omega)={\rm Inv}N(\omega)$, then we have
\begin{equation}\label{2}
\varphi(\omega,S(\omega))=S(\theta\omega)\subset{\rm
int}N(\theta\omega)
\end{equation}
by the fact that $N(\omega)$ is a random isolating block. Since
$\varphi(\omega,\cdot)$ is homeomorphism, by (\ref{1}) and
(\ref{2}), we know that $N^-(\omega)$ and $S(\omega)$ are disjoint
random compact sets. Hence when $L(\omega)$ is a sufficiently small
random compact neighborhood of $N^-(\omega)$, we have
$N(\omega)\backslash L(\omega)$ is a random neighborhood of
$S(\omega)$. By the definitions of random isolating neighborhood and
random isolated invariant set, we know that any random neighborhood
$V(\omega)$ of $S(\omega)$ in $N(\omega)$ is a random isolating
neighborhood with $S(\omega)={\rm Inv}V(\omega)$. In particular,
${\rm cl}(N(\omega)\backslash L(\omega))$ is a random isolating
neighborhood of $S(\omega)$.

Now we verify that
\[
\varphi(\omega,L(\omega))\cap{\rm cl}(N(\theta\omega)\backslash
    L(\theta\omega))=\emptyset.
\]
By Remark \ref{r1} we know that
\[
\varphi(\omega,N(\omega)\backslash N^-(\omega))\subset {\rm
int}N(\theta\omega),
\]
i.e.
\[
N(\omega)\backslash N^-(\omega)\subset\varphi^{-1}(\theta\omega,{\rm
int}N(\theta\omega))\cap N(\omega)\subset
\varphi^{-1}(\theta\omega,N(\theta\omega))\cap N(\omega).
\]
By the closeness of $\varphi^{-1}(\theta\omega,N(\theta\omega))\cap
N(\omega)$, we have
\begin{equation}\label{3}
{\rm cl}(N(\omega)\backslash N^-(\omega))\subset
\varphi^{-1}(\theta\omega,N(\theta\omega))\cap N(\omega).
\end{equation}
On the other hand, by (\ref{N-}) we have
\begin{align}\label{4}
\varphi(\theta_{-1}\omega,N^-(\theta_{-1}\omega))
&=\varphi[\theta_{-1}\omega,
N(\theta_{-1}\omega)\cap(\varphi^{-1}(\omega,{\rm
int}^cN(\omega)))]\nonumber\\
& \subset [\varphi(\theta_{-1}\omega,
N(\theta_{-1}\omega))]\cap[\varphi(\theta_{-1}\omega,\varphi^{-1}(\omega,{\rm
int}^cN(\omega)))]\nonumber\\
& =\varphi(\theta_{-1}\omega, N(\theta_{-1}\omega))\cap {\rm
int}^cN(\omega)
\end{align}
By (\ref{3}) and (\ref{4}) we have
\begin{align*}
&\quad~\varphi(\theta_{-1}\omega,N^-(\theta_{-1}\omega))\cap {\rm
cl}(N(\omega)\backslash N^-(\omega))\\
& \subset\varphi(\theta_{-1}\omega, N(\theta_{-1}\omega))\cap {\rm
int}^cN(\omega)\cap\varphi^{-1}(\theta\omega,N(\theta\omega))\cap
N(\omega)\\
& =\emptyset,
\end{align*}
where `=' holds by the fact that $N(\omega)$ is a random isolating
block. By the fact that $\varphi(\omega,\cdot)$ is a homeomorphism
we have $\varphi(\theta_{-1}\omega,N^-(\theta_{-1}\omega))$ is a
random compact set. Then when $L(\omega)$ is a sufficiently small
random compact neighborhood of $N^-(\omega)$ in $N(\omega)$ we have
\begin{align*}
&\quad~\varphi(\theta_{-1}\omega,L(\theta_{-1}\omega))\cap {\rm
cl}(N(\omega)\backslash L(\omega))\\
& \subset\varphi(\theta_{-1}\omega,L(\theta_{-1}\omega))\cap {\rm
cl}(N(\omega)\backslash N^-(\omega))\\
& =\emptyset.
\end{align*}
Up to now we have proved that $(N(\omega),L(\omega))$ is a random
filtration pair for ${\rm Inv}N(\omega)$.

(2) Assume $\psi$ is a random homeomorphism sufficiently close to
$\varphi$ in the random $C^0$ topology. Then it is easy to see that
\[
\psi(\theta_{-1}\omega, N(\theta_{-1}\omega))\cap N(\omega)\cap
\psi^{-1}(\theta\omega,N(\theta\omega))\subset{\rm int}N(\omega)
\]
by the fact that $N(\omega)$ is a random isolating neighborhood for
$\varphi$. Hence $N(\omega)$ is also a random isolating block for
$\psi$. We denote $N_{\psi}^{-}(\omega)$ the exit set of $N(\omega)$
for $\psi$. In fact, it is easy to verify that $L(\omega)$ being a
random neighborhood of $N^-(\omega)$ in $N(\omega)$ is equivalent to
that
\begin{equation}\label{5}
\varphi(\omega,{\rm cl}(N(\omega)\backslash L(\omega)))\subset {\rm
int}N(\theta\omega)
\end{equation}
by the fact $\varphi(\omega,\cdot)$ is homeomorphism again. Then we
have (\ref{5}) holds if we replace $\varphi$ with $\psi$  when
$\psi$ is sufficiently close to $\varphi$, which indicates that
$L(\omega)$ is a random neighborhood of $N_{\psi}^{-}(\omega)$ in
$N(\omega)$. Therefore $(N(\omega),L(\omega))$ is a random
filtration pair for the random homeomorphisms sufficiently close to
$\varphi$ in the random $C^0$ topology. \hfill $\Box$

Assume $P=(N(\omega),L(\omega))$ is a random filtration pair for
$\varphi$ and denote $N_L(\omega)$ the quotient space
$N(\omega)/L(\omega)$, where the collapsed set $L(\omega)$ is
denoted by $[L(\omega)]$ and is taken as the base-point. That is,
\[
N_L(\omega)=(N(\omega)\backslash
L(\omega)\cup[L(\omega)],[L(\omega)]).
\]
For $\forall\omega\in\Omega$, a set $U(\omega)\subset N_L(\omega)$
is open if either $U(\omega)$ is open in $N(\omega)$ and
$U(\omega)\cap L(\omega)=\emptyset$ or the set
$(U(\omega)\cap(N(\omega)\backslash L(\omega)))\cup L(\omega)$ is
open in $N(\omega)$. In the case $L(\omega)=\emptyset$,
$N_L(\omega)$ is obtained from $N(\omega)$ by adding an isolated
base-point, i.e. $N_L(\omega)=N(\omega)\cup [\emptyset]$. We call
the random quotient space obtained in this way {\em random pointed
space}. In particular, $\emptyset/\emptyset$ is a random pointed
space consisting of just one point--- base point. In the following
we will identify $N_L(\omega)\backslash[L(\omega)]$ with
$N(\omega)\backslash L(\omega)$.

We define the {\em random pointed space map associated to $P$}, i.e.
$\varphi_P(\omega,\cdot):N_L(\omega)\rightarrow N_L(\theta\omega)$,
as follows:
\[
\varphi_P(\omega,\cdot)=\left\{\begin{array}{ll}
                                 [L(\theta\omega)], & x=[L(\omega)]~{\rm or}~
                                 \varphi(\omega,x)\not\in N(\theta\omega),\\
                                 p(\theta\omega,\varphi(\omega,x)), & {\rm
                                 otherwise},
                                \end{array}
\right.
\]
where $p(\omega,\cdot):N(\omega)\rightarrow N_L(\omega)$  is the
quotient map. Then we have the following theorem.

\begin{theorem}\label{th2}
The random pointed space map associated to $P$,
$\varphi_P(\omega,\cdot):N_L(\omega)\rightarrow N_L(\theta\omega)$
is a base-point preserving random map with the property that
$[L(\omega)]\subset {\rm
int}\varphi_P^{-1}(\theta\omega,[L(\theta\omega)])$. Moreover,
$\varphi_P(\omega,\cdot):N_L(\omega)\rightarrow N_L(\theta\omega)$
is continuous, and $\varphi_P(\cdot,x)$ is measurable.
\end{theorem}
{\it Proof.} By the definition of $\varphi_P$ it is obvious that
$\varphi_P$ preserves the base-point. Since
$P=(N(\omega),L(\omega))$ is a random filtration pair, we have
\[
\varphi(\omega,L(\omega))\cap{\rm cl}(N(\theta\omega)\backslash
L(\theta\omega))=\emptyset.
\]
Hence by the compactness of $L(\omega)$ and ${\rm
cl}(N(\theta\omega)\backslash L(\theta\omega))$, there exists a
random neighborhood $K(\omega)$ of $L(\omega)$ such that
\[
\varphi(\omega,K(\omega))\cap{\rm cl}(N(\theta\omega)\backslash
L(\theta\omega))=\emptyset.
\]
Therefore by the definition of $\varphi_P$ we have
$\varphi_P(\omega,x)=[L(\theta\omega)]$ whenever $x\in K(\omega)$,
i.e. $[L(\omega)]\subset {\rm
int}\varphi_P^{-1}(\theta\omega,[L(\theta\omega)])$.

It is obvious that when $x$ satisfies $\varphi(\omega,x)\not\in
N(\theta\omega)$ there exists a small neighborhood $N_x(\omega)$ of
$x$ such that $\varphi(\omega,N_x(\omega))\cap
N(\theta\omega)=\emptyset$. That is, $\varphi_P(\omega,\cdot)$ is
continuous at $x$. When
$\varphi_P(\omega,x)=p(\theta\omega,\varphi(\omega,x))$,
$\varphi_P(\omega,\cdot)$ is the composition of two continuous
functions $p(\theta\omega,\cdot)$ and $\varphi(\omega,\cdot)$, so it
is continuous. Hence we only need to verify that
$\varphi(\omega,\cdot)$ is continuous at $[L(\omega)]$. Assume $x_n$
is a sequence converging to $[L(\omega)]$, then when $n$ is
sufficiently large we have $x\in K(\omega)$. Hence we have
\[
\varphi(\omega,x_n)=[L(\theta\omega)]=\varphi(\omega,[L(\omega)])
\]
when $n$ is sufficiently large, which verifies that
$\varphi(\omega,\cdot)$ is continuous at $[L(\omega)]$.

Since $\varphi(\cdot,x)$ is measurable, we immediately obtain that
$\varphi_P(\cdot,x)$ is measurable by its definition. \hfill $\Box$

\begin{remark}\rm
Similar to deterministic case, to study the dynamical behaviors of
dynamical systems, the definitions of forward invariant, backward
invariant and invariant sets are introduced for RDS, see \cite{Ar1}
for details. Based on this, for a given random set we can define the
exit set of it, the maximal invariant random set in it etc. These
definitions are very natural and have their obvious dynamical
explanations. If we go one step further, inspired by \cite{Fra}, we
can introduce the definition of random filtration pair, which is
crucial to define Conley index for RDS. Of course the random
filtration pair for an isolated invariant set is important to study
the dynamical behaviors for RDS.
\end{remark}

\section{Random shift equivalence}

In order to define the random Conley index for discrete random
dynamical systems, we need to find an invariant for random isolated
invariant sets--random shift equivalence.

Assume $C(\omega),D(\omega)$ are two random pointed spaces and
\[
c(\omega,\cdot):C(\omega)\rightarrow C(\theta\omega),\quad
d(\omega,\cdot):D(\omega)\rightarrow D(\theta\omega)
\]
are two base-point preserving maps satisfying that
$c(\omega,\cdot),d(\omega,\cdot)$ are continuous and
$c(\cdot,x),d(\cdot,x)$ are measurable. If there exist
$r(\omega,\cdot):C(\omega)\rightarrow D(\theta_{n_1}\omega)$,
$s(\omega,\cdot):D(\omega)\rightarrow C(\theta_{n_2}\omega)$, where
$r,s$ have the same properties as $c,d$ (i.e. they satisfy the
continuous-measurable condition and they preserve base-point) and
$n_i=n_i(\omega),i=1,2$ are measurable, such that the following
diagrams
\begin{equation}\label{CD}
\begin{CD}
C(\omega)
@>c(\omega,\cdot)>> C(\theta\omega) \\
@Vr(\omega,\cdot)VV @V({\rm I})Vr(\theta\omega,\cdot)V\\
D(\theta_{n_1}\omega) @>({\rm II})>d(\theta_{n_1}\omega,\cdot)>
D(\theta_*\omega)
\end{CD}
\quad\qquad
\begin{CD}
D(\omega)
@>d(\omega,\cdot)>> D(\theta\omega) \\
@Vs(\omega,\cdot)VV @V({\rm III})Vs(\theta\omega,\cdot)V\\
C(\theta_{n_2}\omega) @>({\rm IV})>c(\theta_{n_2}\omega,\cdot)>
C(\theta_*\omega)
\end{CD}
\end{equation}
are {\em quasi-commutative} and the following holds:
\begin{align}
&r(\theta_{n_2}\omega,s(\omega,\cdot))=d^{n_2(\omega)+n_1(\theta_{n_2}\omega)}(\omega,\cdot),\label{r}\\
&s(\theta_{n_1}\omega,r(\omega,\cdot))=c^{n_1(\omega)+n_2(\theta_{n_1}\omega)}(\omega,\cdot).\label{s}
\end{align}
Then we call $(C,c)$ is {\em random shift equivalent} to $(D,d)$ and
denote it by $(C,c)\sim (D,d)$. We use ``quasi-commutative" meaning
that the two diagrams in (\ref{CD}) are not strictly commutative.
For example, in the first diagram of (\ref{CD}) we obviously have
\begin{align*}
&d(\theta_{n_1}\omega,r(\omega,\cdot)): C(\omega)\rightarrow
D(\theta_{n_1(\omega)+1}\omega),\\
&r(\theta\omega,c(\omega,\cdot)): C(\omega)\rightarrow
D(\theta_{n_1(\theta\omega)+1}\omega).
\end{align*}
That is, (I) and (II) do not generally arrive at the same
destination unless we have $n_1(\omega)=n_1(\theta\omega)$. But in
the situation we confront, this condition is not satisfied
generally. Hence the diagram does not commute generally. But if we
can adjust it to make the diagram commute in the following sense,
i.e. if we have
\begin{equation}\label{co}
\left\{
  \begin{array}{ll}
    r(\theta\omega,c(\omega,\cdot))=d^{n_1(\theta\omega)-n_1(\omega)}
(\theta_{n_1(\omega)+1}\omega,d(\theta_{n_1(\omega)}\omega,
r(\omega,\cdot))), & n_1(\theta\omega)\ge n_1(\omega), \\
  d^{n_1(\omega)-n_1(\theta\omega)}
(\theta_{n_1(\theta\omega)+1}\omega,r(\theta\omega,c(\omega,\cdot)))
=d(\theta_{n_1(\omega)},r(\omega,\cdot)), & n_1(\theta\omega)<
n_1(\omega),
  \end{array}
\right.
\end{equation}
then we call the diagram quasi-commutative. For simplicity, we can
formally rewrite (\ref{co}) as
\[
r(\theta\omega,c(\omega,\cdot))=d^\triangle(\theta_{n_1(\omega)+1}\omega,
d(\theta_{n_1(\omega)}\omega,r(\omega,\cdot))),
\]
where $d^\triangle$ denotes the adjustment we make and it should be
understood as in (\ref{co}), i.e. $d^\triangle$ should vary
according to the relation $n_1(\theta\omega)\ge n_1(\omega)$ or
$n_1(\theta\omega)< n_1(\omega)$. Symmetrically we also can rewrite
(\ref{co}) as
\[
d^\triangle(\theta_{n_1(\theta\omega)+1}\omega,r(\theta\omega,c(\omega,\cdot)))
=d(\theta_{n_1(\omega)}\omega,r(\omega,\cdot)).
\]
And in the first diagram of (\ref{CD}), $D(\theta_*\omega)$ means
that through (I), (II) we may arrive at different destinations, i.e.
$D(\theta_{n_1(\theta\omega)+1}\omega)$ and
$D(\theta_{n_1(\omega)+1}\omega)$. By the above adjustment we write
$D(\theta_*\omega)$ to indicate that the diagram is
quasi-commutative and $\theta_*\omega$ should be valued
appropriately, i.e. $*=\max\{n_1(\omega)+1,n_1(\theta\omega)+1\}$.
Loosely speaking, starting from the same start point ($C(\omega)$),
we can (by adjustment) arrive at the same end point
($D(\theta_*\omega)$) through two paths ($r\circ c$ and $d\circ r$).
This is just our intuitive idea of quasi-commutation. In the
following for simplicity we often use ``$*$" to denote the number we
can obviously understand but we have no need to formulate it
explicitly. The second diagram in (\ref{CD}) can be defined
quasi-commutative in the completely same way. In (\ref{r}), we only
care that the composition of $r$ and $s$ equals the power of $d$,
but we do not care what on earth the power is in spite that we have
explicitly given it. Hence for simplicity, we often formally rewrite
(\ref{r}) as $r\circ s=d^*$ and rewrite (\ref{s}) as $s\circ r=c^*$.

We now verify that the random shift equivalence is an equivalence
relation. The reflexivity and symmetry obviously hold and we only
need to verify the transitivity of random shift equivalence. Assume
$(C,c)$ is random shift equivalent to $(D,d)$ with (\ref{CD}),
(\ref{r}) and (\ref{s}) hold and assume $(D,d)$ is random shift
equivalent to $(E,e)$ with
\begin{align*}
&r_1(\theta\omega,d(\omega,\cdot))=e^\triangle(\theta_{n'_1(\omega)+1}\omega,
e(\theta_{n'_1(\omega)}\omega,r_1(\omega,\cdot))),\\
&s_1(\theta\omega,e(\omega,\cdot))=d^\triangle(\theta_{n'_2(\omega)+1}\omega,
d(\theta_{n'_2(\omega)}\omega,s_1(\omega,\cdot))),\\
&r_1(\theta_{n'_2}\omega,s_1(\omega,\cdot))=e^{n'_2(\omega)+n'_1(\theta_{n'_2}\omega)}(\omega,\cdot),\\
&s_1(\theta_{n'_1}\omega,r_1(\omega,\cdot))=d^{n'_1(\omega)+n'_2(\theta_{n'_1}\omega)}(\omega,\cdot),
\end{align*}
where
\[
r_1(\omega,\cdot):D(\omega):\rightarrow E(\theta_{n'_1}\omega),\quad
s_1(\omega,\cdot):E(\omega):\rightarrow D(\theta_{n'_2}\omega).
\]
If we denote $r_2=r_1\circ d^\triangle\circ r$, then we have
\[
\begin{CD}
C(\omega)@>d^\triangle\circ r\circ c>> D(\theta_{*}\omega) @>r_1>>
E(\theta_{[n'_1(\theta_{*}\omega)+*]}\omega),\\
C(\omega)@>d^\triangle\circ r>> D(\theta_{*}\omega) @>e\circ r_1>>
E(\theta_{[n'_1(\theta_{*}\omega)+*+1]}\omega),
\end{CD}
\]
that is
\begin{align*}
&r_2\circ c: C(\omega)\rightarrow E(\theta_{[n'_1(\theta_{*}\omega)+*]}\omega),\\
&e\circ r_2: C(\omega)\rightarrow
E(\theta_{[n'_1(\theta_{*}\omega)+*+1]}\omega).
\end{align*}
By the assumption $(C,c)\sim (D,d)$ and $(D,d)\sim (E,e)$ it is easy
to see that the following diagram
\begin{equation*}
\begin{CD}
C(\omega)
@>c(\omega,\cdot)>> C(\theta\omega) \\
@Vr_2(\omega,\cdot)VV @VVr_2(\theta\omega,\cdot)V\\
E(\theta_{n^{''}_1}\omega) @>>e(\theta_{n^{''}_1}\omega,\cdot)>
E(\theta_*\omega)
\end{CD}
\end{equation*}
is quasi-commutative, where $n^{''}_1=n'_1(\theta_{*}\omega)+*$ and
$E(\theta_*\omega)$ should be understood similar to
$D(\theta_*\omega)$ above. Let $s_2=s\circ d^\triangle\circ s_1$, in
the completely same way we can obtain that
\[
s_2(\theta\omega,e(\omega,\cdot))=c^\triangle(\theta_{n^{''}_2(\omega)+1}\omega,
c(\theta_{n^{''}_2(\omega)}\omega,s_2(\omega,\cdot))),
\]
that is, the diagram
\begin{equation*}
\begin{CD}
E(\omega)
@>e(\omega,\cdot)>> E(\theta\omega) \\
@Vs_2(\omega,\cdot)VV @VVs_2(\theta\omega,\cdot)V\\
C(\theta_{n^{''}_2}\omega) @>>c(\theta_{n^{''}_2}\omega,\cdot)>
C(\theta_*\omega)
\end{CD}
\end{equation*}
is quasi-commutative. Now we verify that $r_2\circ s_2=e^*$. Notice
that
\begin{align*}
r_2\circ s_2(\omega,\cdot)&=r_1\circ d^\triangle\circ (r\circ
s)\circ d^\triangle\circ
s_1(\omega,\cdot)=r_1\circ d^*\circ s_1(\omega,\cdot)\\
&=(r_1\circ d)\circ d^{*-1}\circ s_1(\omega,\cdot)=(e^\triangle\circ
e\circ
r_1)\circ d^{*-1}\circ s_1(\omega,\cdot)\\
&=e^*\circ(r_1\circ d)\circ d^{*-2}\circ s_1(\omega,\cdot)=\cdots\\
&=e^*\circ r_1\circ s_1(\omega,\cdot)=e^*(\omega,\cdot),
\end{align*}
where ``$*$" is not an invariant number. It is only a number we
obviously understand, and in our eyes we do not distinguish $*$,
$*+\triangle$, $\triangle+1$, $*+*$  etc and hence we denote them by
unified notation ``$*$". In the same way we can obtain $s_2\circ
r_2=c^*$. Hence we have verified $(C,c)\sim (E,e)$, i.e. random
shift equivalence is an equivalence relation.

\begin{lemma}\label{lem3}
Assume $P'=(N(\omega),L'(\omega))$ and $P=(N(\omega)\cup
L(\omega),L(\omega))$ are two random filtration pairs for
$S(\omega)$ and that $L'(\omega)\subset L(\omega)$ and
$\varphi(\omega,L(\omega))\subset {\rm int}L(\theta\omega)$. Then
the induced random maps, $\varphi_{P'}$ and $\varphi_P$ are random
shift equivalent.
\end{lemma}
{\it Proof.} Denote $Q(\omega)=N(\omega)\cup L(\omega)$ and define
$r(\omega,\cdot):N_{L'}(\omega)\rightarrow Q_L(\omega)$ as follows:
\[
r(\omega,\cdot)=\left\{\begin{array}{ll}
                                 [L(\omega)], & x=[L'(\omega)],\\
                                 p(\omega,x), & {\rm
                                 otherwise},
                                \end{array}
\right.
\]
where $p(\omega,\cdot):Q(\omega)\rightarrow Q_L(\omega)$ is the
quotient map. It is easy to see that $r(\omega,\cdot)$ is
continuous, $r(\cdot,x)$ is measurable and
$r(\theta\omega,\varphi_{P'}(\omega,\cdot))=\varphi_P(\omega,r(\omega,\cdot))$.

We affirm that for almost all $\omega\in\Omega$, $\exists
n=n(\omega)$ such that $\varphi^k(\omega,N(\omega)\cap
L(\omega))\subset L'(\theta_k\omega)$ for some $k<n$. (In fact,
similar to Lemma \ref{lem1}, $n(\omega)$ can be chosen measurable.)
If the assertion is false, we have
\begin{align*}
\mathbb P\{\omega|~\exists x\in N(\omega)\cap L(\omega)~{\rm such~
that~}\varphi^k(\omega,x)\in (N(\theta_k\omega)\backslash
L'(\theta_k\omega)),\forall k\in\mathbb N\}>0
\end{align*}
by the fact $L'(\omega)$ is a random neighborhood of $N^-(\omega)$
in $N(\omega)$. Denote this set $\tilde\Omega$, so $\mathbb
P(\tilde\Omega)>0$. For $\forall\omega\in\tilde\Omega$, let
\[
\tilde L(\omega):=(\bigcap_{k\ge
0}\varphi^{-k}(\theta_k\omega,N(\theta_k\omega)\backslash
L'(\theta_k\omega)))\cap(N(\omega)\cap L(\omega)).
\]
Then we have $\tilde L(\omega)\neq\emptyset$ whenever
$\omega\in\tilde\Omega$. Let
\[
\tilde N(\omega):=\left\{
                    \begin{array}{ll}
                      \tilde L(\omega), & \omega\in\tilde\Omega, \\
                      S(\omega), &
\omega\in\Omega\backslash\tilde\Omega.
                    \end{array}
                  \right.
\]
By the definition of $\tilde N(\omega)$ we have
\[
\varphi^k(\omega,\tilde N(\omega))\subset
N(\theta_k\omega)\backslash
L'(\theta_k\omega),\forall\omega\in\Omega,\forall k\in\mathbb N.
\]
Hence $\Omega_{\tilde N}(\omega)$ is an invariant random compact set
in ${\rm cl}(N(\omega)\backslash L'(\omega))$. By the fact that
$S(\omega)$ is the maximal invariant random compact set in ${\rm
cl}(N(\omega)\backslash L'(\omega))$ we have $\Omega_{\tilde
N}(\omega)\subset S(\omega)$ almost surely. But on the other hand,
by the fact $\varphi(\omega,L(\omega))\subset {\rm
int}L(\theta\omega)$ and measure preserving of $\theta_n$ we have
\begin{align*}
\mathbb P\{\omega|~\Omega_{\tilde N}(\omega)\subset
L(\omega)\}&\ge\mathbb P\{\omega|~\varphi^k(\theta_{-k}\omega,\tilde
N(\theta_{-k}\omega))\subset L(\omega),\forall k\in\mathbb
N\}\\
&= \mathbb P\{\omega|~\varphi^k(\omega,\tilde N(\omega))\subset
L(\theta_k\omega),\forall k\in\mathbb N\}\\
&\ge\mathbb P\{\omega|~\tilde N(\omega)\subset L(\omega)\}\\
&=\mathbb P(\tilde\Omega)>0,
\end{align*}
a contradiction to the fact $L(\omega)\cap S(\omega)=\emptyset$ for
$\forall\omega\in\Omega$.

For the above $n(\omega)$, if $n(\omega)\ge n(\theta\omega)$, we do
not make adjustment; if $n(\omega)< n(\theta\omega)$, let
$n(\omega)=n(\theta\omega)$. By such adjustment, $n(\omega)$ is
still measurable and it satisfies $n(\omega)\ge n(\theta\omega)$ and
$\varphi^k(\omega,N(\omega)\cap L(\omega))\subset
L'(\theta_k\omega)$ for some $k<n(\omega)$, and we still denote it
by $n=n(\omega)$.

Noting that $Q(\omega)\backslash L(\omega)=N(\omega)\backslash
L(\omega)$, we define $s(\omega,\cdot):Q_{L}(\omega)\rightarrow
N_{L'}(\theta_n\omega)$ as follows:
\[
s(\omega,\cdot)=\left\{\begin{array}{ll}
                                 [L'(\theta_n\omega)], & x=[L(\omega)],\\
                                 \varphi_{P'}^n(\omega,p'(\omega,x)), & {\rm
                                 otherwise},
                                \end{array}
\right.
\]
where $p'(\omega,\cdot):N(\omega)\rightarrow N_{L'}(\omega)$ is the
quotient map. It is easy to see that
$s(\theta\omega,\varphi_{P}(\omega,\cdot))=\varphi_{P'}^\triangle(\theta_{n+1}\omega,
\varphi_{P'}(\theta_n\omega,s(\omega,\cdot)))$.

It is obvious that $s(\omega,\cdot)$ is continuous on
$Q_L(\omega)\backslash[L(\omega)]$ since it is the composition of
continuous functions. By the fact $\varphi(\omega,\cdot)$ is
homeomorphism there exists a random neighborhood $V(\omega)$ of
$L(\omega)$ such that $\varphi(\omega,V(\omega))\subset{\rm
int}L(\theta\omega)$. Hence for $\forall x\in V(\omega)$, by the
fact $n(\omega)\ge n(\theta\omega)$ we have $\varphi^k(\omega,x)\in
L'(\theta_k\omega)$ for some $k\le n$ and therefore
$\varphi_{P'}^n(\omega,p'(\omega,x))=[L'(\theta_n\omega)]$. Thus we
have obtained $s(\omega,x)$ is continuous at $[L(\omega)]$. The
measurability of $s(\cdot,x)$ is obvious.

By the definitions of $r,s$, it is easy to verify that
\[
s(\omega,r(\omega,\cdot))=\varphi_{P'}^n(\omega,\cdot),
~r(\theta_n\omega,s(\omega,\cdot))=\varphi_P^n(\omega,\cdot).
\]
This completes the proof of the lemma. \hfill $\Box$

\begin{lemma}\label{lem4}
Assume $P=(N(\omega),L(\omega))$ and $P'=(N'(\omega),L(\omega))$ are
two random filtration pairs for $S(\omega)$ and that
$N(\omega)\subset N'(\omega)$ and $\varphi(\omega,L(\omega))\subset
{\rm int}L(\theta\omega)$. Then the induced random maps, $\varphi_P$
and $\varphi_{P'}$ are random shift equivalent.
\end{lemma}
{\it Proof.} Define $r(\omega,\cdot):N_{L}(\omega)\rightarrow
N'_L(\omega)$ as follows:
\[
r(\omega,\cdot)=\left\{\begin{array}{ll}
                                 [L(\omega)], & x=[L(\omega)],\\
                                 p'(\omega,x), & {\rm
                                 otherwise},
                                \end{array}
\right.
\]
where $p'(\omega,\cdot):N'(\omega)\rightarrow N'_L(\omega)$ is the
quotient map. It is easy to see that $r(\omega,\cdot)$ is
continuous, $r(\cdot,x)$ is measurable and
$r(\theta\omega,\varphi_{P}(\omega,\cdot))=\varphi_{P'}(\omega,r(\omega,\cdot))$.

Since $L(\omega)$ is a random neighborhood of $N^-(\omega)$, we have
$ \varphi(\omega,N(\omega)\backslash L(\omega))\subset{\rm
int}N(\theta\omega). $ This together with the assumption
$\varphi(\omega,L(\omega))\subset {\rm int}L(\theta\omega)$ imply
that $\varphi(\omega,N(\omega))\subset{\rm int}N(\theta\omega)$,
i.e. $N(\omega)$ is forward invariant. Hence we have
$\Omega_{N}(\omega)\subset{\rm int}N(\omega)$ and it is the maximal
invariant random compact set in $N(\omega)$.

We can decompose the random set $N'(\omega)\backslash N(\omega)$
into two random sets: $N_1(\omega)$ and $N_2(\omega)$. They are
determined as follows:
\begin{align*}
& N_1(\omega):=\{x|~x\in N'(\omega)\backslash
N(\omega),\varphi^k(\omega,x)\in N'(\theta_k\omega)\backslash
L(\theta_k\omega),\hbox{for}~\forall k\in\mathbb
N\},\\
& N_2(\omega):=\{x|~x\in N'(\omega)\backslash N(\omega),\exists
k\in\mathbb N ~\hbox{such that}~ \varphi^k(\omega,x)\in
L(\theta_k\omega)\}.
\end{align*}
By the definition of $N_1(\omega)$ and the measure preserving of
$\theta_n$ we obtain that $ \Omega_{N_1}(\omega)\subset{\rm
cl}(N'(\omega)\backslash L(\omega)) $. By the definition of
$N_2(\omega)$ and the forward invariance of $L(\omega)$ we have
$\Omega_{N_2}(\omega)\subset L(\omega)$. Since $S(\omega)$,
$\Omega_L(\omega)$ are maximal invariant random set in ${\rm
cl}(N'(\omega)\backslash L(\omega))$, $L(\omega)$ respectively, we
obtain that
\[
\Omega_{N'\backslash
N}(\omega)=\Omega_{N_1}(\omega)\cup\Omega_{N_2}(\omega)\subset
S(\omega)\cup\Omega_L(\omega)\subset{\rm int}N(\omega).
\]
Therefore,
\[
\Omega_{N'}(\omega)=\Omega_{N}(\omega)\cup\Omega_{N'\backslash
N}(\omega)\subset{\rm int}N(\omega).
\]
Hence for $\forall\omega\in\Omega$, $\exists n(\omega)$ ( in fact
$n(\omega)$ may be chosen measurable similar to Lemma \ref{lem1})
such that
\[
\varphi^n(\theta_{-n}\omega,N'(\theta_{-n}\omega))\subset {\rm
int}N(\omega)~{\rm whenever~}n\ge n(\omega)-1
\]
by the forward invariance of $N'(\omega)$ (the forward invariance of
$N'(\omega)$ follows completely similar to the proof of forward
invariance of $N(\omega)$). In fact, it is easy to verify that for
arbitrary forward invariant random set $D(\omega)$, we have
\[
\varphi^n(\theta_{-n}\omega,D(\theta_{-n}\omega))\subset
\varphi^m(\theta_{-m}\omega,D(\theta_{-m}\omega)),~\hbox{whenever}~n\ge
m.
\]
Then by the measure preserving of $\theta_n$ we have
\begin{align*}
&\mathbb P\{\omega|~\varphi^n(\omega,N'(\omega))\subset {\rm
int}N(\theta_n\omega),n\ge n(\omega)\}\\
=&\mathbb
P\{\omega|~\varphi^n(\theta_{-n}\omega,N'(\theta_{-n}\omega))\subset
{\rm int}N(\omega),n\ge n(\omega)\}\\
=&1.
\end{align*}

Adjust the above $n(\omega)$ similar to in Lemma \ref{lem3} if
necessary such that $n(\omega)\ge n(\theta\omega)$ and still denote
it by $n=n(\omega)$. Define
$s(\omega,\cdot):N'_{L}(\omega)\rightarrow N_{L}(\theta_n\omega)$ as
follows:
\[
s(\omega,\cdot)=\left\{\begin{array}{ll}
                                 [L(\theta_n\omega)], & x=[L(\omega)],\\
                                 p(\theta_n\omega,\varphi^n(\omega,x)), & {\rm
                                 otherwise},
                                \end{array}
\right.
\]
where $p(\omega,\cdot):N(\omega)\rightarrow N_{L}(\omega)$ is the
quotient map. It is easy to see that
$s(\theta\omega,\varphi_{P'}(\omega,\cdot))=\varphi_{P}^\triangle(\theta_{n+1}\omega,
\varphi_{P}(\theta_n\omega,s(\omega,\cdot)))$. It is obvious that
$s(\omega,\cdot)$ is continuous on
$N'_L(\omega)\backslash[L(\omega)]$ since it is the composition of
continuous functions. By the assumption
$\varphi(\omega,L(\omega))\subset {\rm int}L(\theta\omega)$ there
exists a random neighborhood $V(\omega)$ of $L(\omega)$ such that
$\varphi(\omega,V(\omega))\subset{\rm int}L(\theta\omega)$. Then by
the fact $n(\omega)\ge n(\theta\omega)$ we have
$s(\omega,V(\omega))=[L(\theta_n\omega)]$, hence $s(\omega,\cdot)$
is continuous at $[L(\omega)]$. The measurability of $s(\cdot,x)$ is
obvious.

By the definitions of $r,s$, it is easy to verify that
\[
r(\theta_n\omega,s(\omega,\cdot))=\varphi_{P'}^n(\omega,\cdot),~
s(\omega,r(\omega,\cdot))=\varphi_{P}^n(\omega,\cdot).
\]
This completes the proof of the lemma. \hfill $\Box$

By the above two lemmas we can now show that the random shift
equivalent class of random pointed space maps is an invariant for a
given random isolated invariant set, which is of crucial importance
for the definition of random Conley index.

\begin{theorem}\label{se}
Assume $P=(N(\omega),L(\omega))$ and $P'=(N'(\omega),L'(\omega))$
are two random filtration pairs for $S(\omega)$, then the induced
random maps, $\varphi_P$ and $\varphi_{P'}$ on the corresponding
random pointed spaces, are random shift equivalent.
\end{theorem}
{\it Proof.} By Lemmas \ref{lem1} and \ref{lem2} we can choose
$\epsilon(\omega)$ sufficiently small such that
\[
{\rm cl}(C_\epsilon(N\backslash L,S)(\omega))={\rm
cl}(C_\epsilon(N'\backslash L',S)(\omega))\subset{\rm
int}(N(\omega)\backslash L(\omega))\cap{\rm
int}(N'(\omega)\backslash L'(\omega)).
\]
Denote $B(\omega):={\rm cl}(C_\epsilon(N\backslash L,S)(\omega))$
and assume $B_0(\omega)$ is a sufficiently small random neighborhood
of $B^-(\omega)$ in $B(\omega)$, then $P_0=(B(\omega),B_0(\omega))$
is random filtration pair for $S(\omega)$ by Theorem \ref{th1}. We
will verify that $\varphi_{P_0}$ is random shift equivalent to
$\varphi_{P}$ and $\varphi_{P'}$. Since the proof is the same we
only need to prove $\varphi_{P_0}$ is random shift equivalent to
$\varphi_{P}$.

We first prove that for $\forall\omega\in\Omega$, $\exists
n=n(\omega)$ such that
\begin{equation}\label{B}
\varphi_P^n(\omega,B^-(\omega))=[L(\theta_n\omega)].
\end{equation}
Since $B(\omega)\subset N(\omega)\backslash L(\omega)$ and we
identify $N_L(\omega)\backslash[L(\omega)]$ with
$N(\omega)\backslash L(\omega)$, we may consider $B(\omega)$ as a
subset of $N_L(\omega)\backslash[L(\omega)]$ and consider $\varphi$
to be $\varphi_P$. By the definition of $B^-(\omega)$, for arbitrary
random variable $x(\omega)\in B^-(\omega)$, there is no random
$\epsilon(\omega)$-chain from
$\varphi_P(\theta_{-1}\omega,x(\theta_{-1}\omega))$ to $S(\omega)$
with positive probability. Since if there is one such
$\epsilon(\omega)$-chain, then we obtain that
$\varphi_P(\theta_{-1}\omega,x(\theta_{-1}\omega))\subset
C_\epsilon(N\backslash L,S)(\omega)$ with positive probability, a
contradiction to the definition of $B^-(\omega)$. Hence we have
$\Omega_{B^-}(\omega)\cap S(\omega)=\emptyset$ almost surely. By the
fact $\varphi_P(\omega,\cdot):N_L(\omega)\rightarrow
N_L(\theta\omega)$ and the compactness of $N_L(\omega)$ we obtain
that $\Omega_{B^-}(\omega)\neq\emptyset$ almost surely and
$\Omega_{B^-}(\omega)\subset N_L(\omega)$. It is easy to see that
there is two invariant random set in $N_L(\omega)$ under the
iteration of $\varphi_P$: $S(\omega)$ and $\{[L(\omega)]\}$. Hence
we obtain that $\Omega_{B^-}(\omega)=\{[L(\omega)]\}$. By Theorem
\ref{th2} we have $[L(\omega)]\subset{\rm
int}\varphi_P^{-1}(\theta\omega,[L(\theta\omega)])$. This together
with the measure preserving of $\theta_n$ imply that for almost all
$\omega\in\Omega$, $\exists n(\omega)$ such that
\begin{equation}\label{11}
\varphi_P^n(\omega,B^-(\omega))\subset\tilde L(\theta_n\omega),~n\ge
n(\omega),
\end{equation}
where $\tilde L(\omega):={\rm
int}\varphi_P^{-1}(\theta\omega,[L(\theta\omega)])$. Hence we have
verified (\ref{B}). Since $\tilde L(\omega)$ is forward invariant
under the iteration of $\varphi_P$, similar to Lemma \ref{lem4},
$n(\omega)$ in (\ref{B}) can be chosen measurable.

Since $B_0(\omega)$ is a sufficiently small random neighborhood of
$B^-(\omega)$, by (\ref{11}) we obtain that
\begin{equation}\label{B0}
\varphi_P^n(\omega,B_0(\omega))=[L(\theta_n\omega)],
\end{equation}
where $n=n(\omega)$ measurable. Adjust this $n(\omega)$ similar to
in Lemmas \ref{lem3} and \ref{lem4} such that $n(\omega)\le
n(\theta\omega)$ and still denote it by $n=n(\omega)$. Let
\[
K(\omega):={\rm cl}({\rm
int}\varphi_P^{-n}(\theta_n\omega,[L(\theta_n\omega)]))\cap
N_L(\omega).
\]
Then by the fact $[L(\omega)]\subset{\rm
int}\varphi_P^{-1}(\theta\omega,[L(\theta\omega)])$ and the fact
$n(\omega)\le n(\theta\omega)$, we obtain that
$\varphi_P(\omega,K(\omega))\subset {\rm int}K(\theta\omega)$.
Moreover, by the definition of $K(\omega)$ we have
$B_0(\omega)\subset{\rm int}K(\omega)$.

Let $Q=(B(\omega)\cup K(\omega),K(\omega))$,
$R=(N_L(\omega),K(\omega))$, then it is easy to see that $Q,R$ are
random filtration pairs for $S(\omega)$. By Lemma \ref{lem3} we
obtain that $\varphi_{P_0}\sim\varphi_{Q}$ and by Lemma \ref{lem4}
we have $\varphi_{Q}\sim\varphi_{R}$. Hence
$\varphi_{P_0}\sim\varphi_{R}$.

Let $\tilde R=(N(\omega),p^{-1}(\omega,K(\omega)))$, where
$p(\omega,\cdot):N(\omega)\rightarrow N_L(\omega)$ is the random
quotient map. Then by Lemma \ref{lem3} we obtain that
$\varphi_{P}\sim\varphi_{\tilde R}$. Noting that we may identify the
random pointed spaces and the corresponding random pointed space
maps of the random filtration pairs $\tilde R$ and $R$, i.e.
$\varphi_{R}\sim\varphi_{\tilde R}$. Hence
$\varphi_{P}\sim\varphi_{R}$. Therefore we have obtained
$\varphi_{P_0}\sim\varphi_{P}$. This completes the proof of the
theorem. \hfill $\Box$

\begin{remark}\rm
To obtain an invariant for a given isolated invariant set, parallel
to \cite{Fra}, we introduced the definition of random shift
equivalence. In fact, it can be regarded as a generalization of
conjugacy of two RDSs. (For the definition of conjugacy of two RDSs,
we refer to \cite{Ar1,IL,IS,LL,LL2} for details and related
applications.) To see this, notice that in the definition of random
shift equivalence, in particular if we have
$r(\omega,\cdot):C(\omega)\rightarrow D(\omega)$ and
$s(\omega,\cdot):D(\omega)\rightarrow C(\omega)$, then $rs={\rm Id}$
and $r$ (or $s$) play the role of cohomological random homeomorphism
between $c,d$, i.e. $c$ is conjugate to $d$ through $r$ (or $s$).

 As pointed out in \cite{Fra}, the random shift
equivalence constructed in Theorem \ref{se} is not unique, here we
only give one relative simple construction. In the definition of
random shift equivalence, for simplicity, we have required that the
random maps $r$ and $s$ preserve base-point, which of course is not
necessary. In other words, the definition of random shift
equivalence can be made slightly more general as that given in
\cite{Fra} for deterministic case.
\end{remark}

\section{Definition of Conley index for RDS}

Assume $S(\omega)$ is a random isolated invariant set and
$P=(N(\omega),L(\omega))$ is a random filtration pair for
$S(\omega)$, by Theorem \ref{se}, we know that the random shift
equivalent class of $\varphi_P:N_L(\omega)\rightarrow
N_L(\theta\omega)$ is an invariant for $S(\omega)$. It is well-known
that continuation of Conley index is one of its most important
properties, hence before defining a random Conley index we need to
explain what random continuation is. Assume $C(\omega)$ is a random
pointed space and $f(\omega,\cdot):C(\omega)\rightarrow
C(\theta_{n}\omega)$, where $n=n(\omega)$ is measurable, is a
base-point preserving random map such that $f(\omega,\cdot)$ is
continuous and $f(\cdot,x)$ is measurable. Let
$g(\omega,\cdot):C(\omega)\rightarrow C(\theta_{n}\omega)$ is
another random map with the same property as that of $f$. We call
$f$ is {\em random homotopic to} $g$, denoted by $f\simeq g$, if
there exist $H:[0,1]\times\Omega\times C(\omega)\rightarrow
C(\theta_n\omega)$ satisfying that $H(\cdot,\omega,\cdot)$ is
continuous and $H(t,\cdot,x)$ is measurable. Moreover, $H$ satisfies
\[
\left\{
  \begin{array}{ll}
     H(0,\cdot,\cdot)=f(\cdot,\cdot),\\
     H(1,\cdot,\cdot)=g(\cdot,\cdot).
  \end{array}
\right.
\]
It is easy to see that random homotopy is an equivalence relation.
We denote $[f]$ the random homotopy class with $f$ the
representative element.

Assume $C(\omega)$ and $D(\omega)$ are two random pointed spaces and
\[
[c(\omega,\cdot)]:C(\omega)\rightarrow C(\theta\omega),\quad
[d(\omega,\cdot)]:D(\omega)\rightarrow D(\theta\omega)
\]
are random homotopy classes. If there exist random homotopy classes
\[
[r(\omega,\cdot)]: C(\omega)\rightarrow D(\theta_{n_1}\omega),\quad
[s(\omega,\cdot)]: D(\omega)\rightarrow C(\theta_{n_2}\omega),
\]
where $n_i=n_i(\omega),i=1,2$ are measurable, such that
\begin{align*}
&r(\theta\omega,c(\omega,\cdot))\simeq
d^\triangle(\theta_{n_1(\omega)+1}\omega,
d(\theta_{n_1(\omega)}\omega,r(\omega,\cdot))),\\
&s(\theta\omega,d(\omega,\cdot))\simeq
c^\triangle(\theta_{n_2(\omega)+1}\omega,
c(\theta_{n_2(\omega)}\omega,s(\omega,\cdot))),\\
&r(\theta_{n_2}\omega,s(\omega,\cdot))\simeq d^{n_2(\omega)+n_1(\theta_{n_2}\omega)}(\omega,\cdot),\\
&s(\theta_{n_1}\omega,r(\omega,\cdot))\simeq
c^{n_1(\omega)+n_2(\theta_{n_1}\omega)}(\omega,\cdot),
\end{align*}
where $d^\triangle$, $c^\triangle$ denote the adjustment we make,
then we call the random homotopy classes $[C,c]$ and $[D,d]$ are
random shift equivalent.

Now we are in the position to give the definition of random Conley
index for random isolated invariant sets.

\begin{definition}\label{Con}\rm
Assume $\varphi$ is the time one map of a discrete random dynamical
system, $S(\omega)$ is a random isolated invariant set for $\varphi$
and $P=(N(\omega),L(\omega))$ is a random filtration pair for
$S(\omega)$. Denote $h_P(S,\varphi)$ the random homotopy class
$[\varphi_P]$ on the random pointed space $N_L(\omega)$ with
$\varphi_P$ a representative element, then we define the {\em random
Conley index $h(S,\varphi)$ for $S(\omega)$} to be the random shift
equivalent class of $h_P(S,\varphi)$.
\end{definition}

\begin{remark}\rm
It is well known that the deterministic Conley index is defined for
autonomous flows or maps, but for discrete RDS we confront, the
random homeomorphism is non-autonomous. As for Conley index for
non-autonomous flows, some authors have dealt with them, see
\cite{War,War1,War2,Pri} among others. They treat them as
skew-product flows by the method of extending phase space. But their
method is inappropriate for our purpose in spite that a discrete RDS
can be seen as a measurable discrete skew-product flow. In fact,
their method depends crucially on the fact that the hull of vector
field function $H(f)$ is compact, see \cite{War,War1,War2,Pri} for
details. But for RDS, the probability space $(\Omega, \mathscr F)$
is only a measurable space and we can not in general assume that it
is a compact topology space, which is too restrictive for
applications.
\end{remark}

The following two theorems state that the random Conley index for
discrete RDS  has the similar properties to that of Conley index for
deterministic maps.

\begin{theorem}{\rm (Continuation property)} Assume $\varphi_\lambda,\lambda\in[0,1]$ is a family of
random homeomorphisms, which depends continuously on $\lambda$. If
$N(\omega)$ is a random isolating neighborhood for each
$\varphi_\lambda,\lambda\in[0,1]$, then the random Conley index
$h(S_\lambda,\varphi_\lambda)$ for $\varphi_\lambda$ is independent
of $\lambda\in[0,1]$, i.e.
$h(S_\lambda,\varphi_\lambda)=h(S_0,\varphi_0)$, where
$S_\lambda(\omega):={\rm Inv}(N(\omega),\varphi_\lambda)$,
$\lambda\in[0,1]$ stands for the random isolated invariant set for
$\varphi_\lambda$ in $N(\omega)$.
\end{theorem}
{\it Proof.} The proof follows immediately from Theorem \ref{th1}.
\hfill $\Box$

We use $\underline 0$ to denote the random Conley index of random
pointed spaces consisting of just one random point with random
constant maps as their corresponding random pointed space maps.
Given a random isolated invariant set $S(\omega)$ in a random
pointed space $N_L(\omega)$ with $\varphi_P$ the corresponding
random pointed space map, by the definition of random shift
equivalence it is easy to verify that the random Conley index
$h(S,\varphi)$ for $S(\omega)$ is $\underline 0$ if and only if
$\varphi_P\simeq f$ and $\varphi_P^n\simeq f^n$ for some random map
$f(\omega,\cdot):N_L(\omega)\rightarrow N_L(\theta\omega)$ and some
measurable $n=n(\omega)$, where $f^n:N_L(\omega)\rightarrow
N_L(\theta_n\omega)$ is the base-point valued random constant map.

The random Conley index can be used to study the structure of random
invariant set, see the following theorem.

\begin{theorem}{\rm (Waz\.ewski property)} Assume $S(\omega)$ is a random isolated
invariant set and the random Conley index for $S(\omega)$ is not
trivial, i.e. $h(S,\varphi)\neq \underline 0$, then
$S(\omega)\neq\emptyset$ almost surely when $\theta$ is ergodic
under $\mathbb P$.
\end{theorem}
{\it Proof.} If $S(\omega)=\emptyset$ with positive probability,
then $S(\omega)=\emptyset$ almost surely by the ergodicity of
$\theta$ and the invariance of $S(\omega)$.  Hence we have
$(\emptyset,\emptyset)$ is a random filtration pair for $S(\omega)$.
It follows that $h(S,\varphi)=\underline 0$, a contradiction. \hfill
$\Box$

\section{Relations between time-continuous RDS and discrete one}

Assume $\phi$ is a time-continuous RDS and consider its time-$h$
map---$\varphi_h(\omega):=\phi(h,\omega):X\rightarrow X$, where
$h>0$. Then $\varphi_h^k(\omega)=\phi(kh,\omega)$ and $\varphi_h$ is
a discrete RDS (recalling that throughout the paper we identify a
random homeomorphism with the discrete RDS generated by it)
generated by the time-$h$ map of the time-continuous RDS $\phi$. For
time-continuous RDS, replacing $n\in\mathbb Z$ by $t\in\mathbb R$,
we can introduce isolated invariant set and isolating neighborhood
similar to Definition \ref{def}.

In this section we simply discuss the relation of isolated invariant
sets between time-continuous RDS and the discrete one generated by
its time-$h$ map.

\begin{theorem}
Assume $\phi$ is a time-continuous RDS and there exists a $\delta>0$
such that for $\forall h\in(0,\delta]$, $S(\omega)$ is an isolated
invariant set of $\varphi_h$, then $S(\omega)$ is an isolated
invariant set of $\phi$.
\end{theorem}
{\it Proof.} It is obvious that $S(\omega)$ is invariant for $\phi$
by the fact that $S(\omega)$ is invariant with respect to
$\varphi_h$ for $\forall h\in(0,\delta]$, so we only need to show it
is isolated with respect to $\phi$. Assume $N(\omega)$ is an
isolating neighborhood of $S(\omega)$ with respect to $\varphi_h$
for some $h\in(0,\delta]$. By the fact
\[
S(\omega)\subset{\rm Inv}(N(\omega),\phi)\subset {\rm
Inv}(N(\omega),\varphi_h)=S(\omega)
\]
we immediately obtain the desired result. \hfill $\Box$

\begin{theorem}
Assume $\phi$ is a time-continuous RDS and $S(\omega)$ is an
isolated invariant set of $\phi$. Moreover, if there exists an
isolating neighborhood $N(\omega)$ of $S(\omega)$ such that the map
\begin{equation}\label{cont}
t\rightarrow N(\theta_t\omega)\in\mathscr K(X)
\end{equation}
is continuous (in general by the definition of random compact set
the map is only measurable), where $\mathscr K(X)$ denotes the space
of non-empty compact subsets of $X$ endowed with the Hausdorff
metric. Then for $\forall h>0$, $S(\omega)$ is an isolated invariant
set of $\varphi_h$.
\end{theorem}
{\it Proof.} The idea of proof is originated from \cite{Mr}. For
fixed $h>0$, by the fact $S(\omega)$ is isolated (with respect to
$\phi$) by $N(\omega)$, the compactness of $[0,h]$ and the
continuity of the map (\ref{cont}) we can choose a compact
neighborhood $\tilde N(\omega)\subset N(\omega)$ of $S(\omega)$ such
that
\[
\phi(t,\omega)\tilde N(\omega)\subset N(\theta_t\omega), ~\forall
t\in[0,h].
\]
In fact we can choose
\[
\tilde
N(\omega)=\bigcap_{t\in[0,h]}\phi(-t,\theta_t\omega)N(\theta_t\omega)
=\bigcap_{t\in[0,h]\cap\mathbb
Q}\phi(-t,\theta_t\omega)N(\theta_t\omega),
\]
where the second ``=" holds by the continuity of the map
(\ref{cont}). Clearly $\tilde N(\omega)$ obtained in this way is a
random compact set. We only need to show that $S(\omega)$ is the
maximal invariant random compact set in $\tilde N(\omega)$ with
respect to $\varphi_h$, i.e.
\[
S(\omega)={\rm Inv}(\tilde N(\omega),\varphi_h).
\]
For $\forall x\in{\rm Inv}(\tilde N(\omega),\varphi_h)$, we have
\[
\varphi^k_h(\omega,x)\in\tilde N(\theta_{kh}\omega), ~\forall
k\in\mathbb Z.
\]
By the choice of $\tilde N(\omega)$ we easily obtain that
\[
\phi(t,\omega,x)\in N(\theta_t\omega),~\forall t\in\mathbb R,
\]
which implies that $x\in S(\omega)$. Hence we obtained the desired
result. \hfill $\Box$

\begin{remark}\rm
Of course assuming the map (\ref{cont}) being continuous is very
restrictive for applications, but there is some cases the assumption
does hold. For instance, when an RDS admits a random attractor (in
the sense of \cite{Sch}) then there is compact neighborhood of the
attractor satisfying this assumption which play the role of
isolating neighborhood, noting that for any invariant random compact
set the continuity of the map (\ref{cont}) naturally holds.
\end{remark}

\section{Several simple examples}
In this section we give several simple examples to illustrate our
results.
\begin{example}\rm
Assume the $2\times 2$ random matrix $A(\omega)=\left(
                    \begin{array}{cc}
                      a(\omega) & 0 \\
0 &b(\omega)     \end{array} \right) $ is the time-one map of a
discrete random dynamical system with phase space $\mathbb R^2$ and
assume $N(\omega)\equiv\{(x,y)|~x^2+y^2\le 1\}$.\\
(i) If $0<a(\omega)<1,0<b(\omega)<1$ for $\forall\omega\in\Omega$,
then it is easy to see (in fact, we have
$A(\omega)N(\omega)\subset{\rm int}N(\theta\omega)$) that
$N(\omega)$ is a random isolating block and $S(\omega)\equiv \{0\}$
is the corresponding random isolated invariant set. It is obvious
that $N^-(\omega)\equiv\emptyset$ and $P=(N(\omega),\emptyset)$ is a
random filtration pair for $S(\omega)$. It is easy to see that any
power of the random pointed space map $A_P$ on the random pointed
space $N_L(\omega)=N(\omega)\cup[\emptyset]$ is not random homotopic
to the base-point valued random constant map, hence we have
$h(S,A)\neq\underline0$.\\
(ii) If $a(\omega)>1,b(\omega)>1$ for $\forall\omega\in\Omega$, then
$N(\omega)$ is a random isolating block and $S(\omega)\equiv \{0\}$
is the corresponding random isolated invariant set. Assume
$P=(N(\omega),L(\omega))$ is a random filtration pair for
$S(\omega)$, where $L(\omega)$ is a sufficiently small neighborhood
of $N^-(\omega)$ (it is obvious that $N^-(\omega)\neq \emptyset$ for
all $\omega\in\Omega$). The random pointed space $N_L(\omega)$ can
be regarded as $S^2$ for each $\omega$ and it is easy to see that
$A_P(\omega):N_L(\omega)\rightarrow N_L(\theta\omega)$ is a
surjection. Hence any power of $A_P$ is not random homotopic to the
base-point valued random constant map, i.e.
$h(S,A)\neq\underline0$.\\
(iii) If $a(\omega)>1,0<b(\omega)<1$ for $\forall\omega\in\Omega$
and assume $N(\omega)$ is a random isolating block for
$S(\omega)\equiv \{0\}$ (it is easy to see that there is possibility
that $N(\omega)$ is not a random isolating block for $S(\omega)$, so
we make such assumption). It is easy to see that
$N^-(\omega)\neq\emptyset$ for each $\omega$ and assume $L(\omega)$
is a sufficiently small neighborhood of $N^-(\omega)$, then
$P=(N(\omega),L(\omega))$ is a random filtration pair for
$S(\omega)$. Then by a simple verification we can obtain that
$h(S,A)\neq\underline0$.
\end{example}

\begin{example}\label{examp}\rm
Consider the logistic model:
\[
\dot{x}=rx(1-\frac{x}{K}),
\]
where $0<r\le 1, K>0$ are constants. By Euler approximation method
we can obtain the corresponding logistic difference equation as
follows:
\[
x_{n+1}=x_n+hrx_n(1-\frac{x_n}{K}),
\]
where $h$ is the step size of Euler method. Now assume that $r$ is
perturbed by a real noise, say,
\[
r(\omega)=r+\xi(\omega)
\]
with $|\xi(\omega)|<r$. We now consider some properties of the
perturbed difference equations:
\[
x_{n+1} =x_n+hr(\theta_{nh}\omega)x_n(1-\frac{x_n}{K}).
\]
 To this end, consider the
following parameterized difference equations:
\begin{equation}\label{exa}
x_{n+1}=x_n+hr_\lambda(\theta_{nh}\omega)x_n(1-\frac{x_n}{K})
\end{equation}
with $r_\lambda(\omega)=r+\lambda\xi(\omega)$, $\lambda\in[0,1]$.
Then for each $\lambda\in[0,1]$, (\ref{exa}) generates a discrete
RDS, denoting it by $\varphi_\lambda$. For the generation of
discrete RDS from random difference equation, the reader can refer
to \cite{Ar1} for details.
Set $X=[0,+\infty)$.\\
(1) Let $N=[0,M]$ with $M>K$. Since $r_\lambda(\omega)>0$ for all
$\omega$, we have
\[
\varphi_\lambda(h,\omega)N\subset{\rm int}N
\]
for small $h$, which implies that $N$ is an isolating neighborhood
for the family $\varphi_\lambda$, $\lambda\in[0,1]$. It is clear
that when $\lambda=0$, i.e. for the unperturbed deterministic
difference equation, the Conley index for the isolated invariant set
$S_0=[0,K]$ is nontrivial. Hence by the continuity property of
random conley index we have that
$h(S_\lambda,\varphi_\lambda)=h(S_0,\varphi_0)\neq\underline 0$,
$\lambda\in[0,1]$. Therefore, by the Waz\.ewski property of random
Conley index we obtain that there exists an invariant random compact
set in $N$ for the perturbed difference equation.\\
(2) Consider any closed interval $\tilde N$ containing $K$  but not
containing $0$. Then it is easy to verify that
\[
\varphi_\lambda(h,\omega)\tilde N\subset{\rm int}\tilde N
\]
when $h$ is small, which implies that $\tilde N$ is an isolating
neighborhood for the family $\varphi_\lambda$, $\lambda\in[0,1]$.
Therefore, by the similar argument to that of (1) we obtain that
there exists an invariant random compact set in $\tilde N$ for the
perturbed difference equation for small $h$. Since the interval
$\tilde N$ can be chosen sufficiently small, we can conclude that
$\{K\}$ is an invariant random compact set for the perturbed
difference equation, which
coincides with our intuition. \\
(3) Consider any closed interval like $[0,a]$ with $a<K$. Denote
$\hat N=[0,a]$, then it is easy to see that
\[
[0,K)\supset\varphi_\lambda(h,\omega){\rm int}\hat N\supset\hat N
\]
for small $h$, which implies that $\hat N$ is an isolating
neighborhood for the family $\varphi_\lambda$, $\lambda\in[0,1]$.
Again by the same argument to that of (1), (2) we obtain that
$\{0\}$ is an invariant random compact set for the perturbed
difference equation for small $h$, which coincides also with our
intuition.
\end{example}

\begin{example}\rm
Consider the Lorenz system in $\mathbb R^3$ described by the
equations:
\[
\left\{
  \begin{array}{l}
    \dot x=\sigma(y-x),\\
    \dot y=\rho x-y-xz, \\
    \dot z=xy-\beta z
  \end{array}
\right.
\]
with $\sigma,\rho,\beta>0$. By Euler approximation method we obtain
the corresponding discrete Lorenz system:
\begin{equation}\label{unp}
X_{n+1}=(hB+I)X_n+hF(X_n),
\end{equation}
where $h$ is the step size of Euler approximation method,
$X_n=(x_n,y_n,z_n)^\top$, $I$ is the identity matrix, $B=\left(
         \begin{array}{ccc}
             -\sigma & \sigma & 0 \\
             \rho & -1 & 0\\
             0 & 0 & -\beta \\
             \end{array}
    \right)
$, and $F(X_n)=(0,-x_nz_n,x_ny_n)^\top$. Assume that
$\sigma,\rho,\beta$ is perturbed by real noises, say,
\begin{align*}
&\sigma(\omega)=\sigma+\xi(\omega),\\
&\rho(\omega)=\rho+\eta(\omega),\\
&\beta(\omega)=\beta+\zeta(\omega).
\end{align*}
Denote $\sigma_\lambda(\omega)=\sigma+\lambda\xi(\omega)$,
$\rho_\lambda(\omega)=\rho+\lambda\eta(\omega)$,
$\beta_\lambda(\omega)=\beta+\lambda\zeta(\omega)$,
$\lambda\in[0,1]$. Consider the family of parameterized difference
systems:
\begin{equation}\label{ex}
X_{n+1}=(hB_\lambda(\theta_{nh}\omega)+I)X_n+hF(X_n)
\end{equation}
with
\[
B_\lambda(\theta_{nh}\omega)=\left(
         \begin{array}{ccc}
             -\sigma_\lambda(\theta_{nh}\omega) & \sigma_\lambda(\theta_{nh}\omega) & 0 \\
             \rho_\lambda(\theta_{nh}\omega) & -1 & 0\\
             0 & 0 & -\beta_\lambda(\theta_{nh}\omega) \\
             \end{array}
    \right).
\]
When $\lambda=0$, (\ref{ex}) corresponds to the unperturbed
discretized Lorenz system by Euler approximation method; when
$\lambda=1$, (\ref{ex}) corresponds to the discrete Lorenz system
perturbed by real noise. Similar to Example \ref{examp}, (\ref{ex})
generates a discrete RDS for each $\lambda$ and denote it by $\varphi_\lambda$.  \\
(1) Assume that for each $\lambda\in[0,1]$,
$0<\rho_\lambda(\omega)<\sigma_\lambda(\omega)\le 1$ almost surely.
Noting that
\begin{align*}
\langle X_{n+1},X_n\rangle
&=\langle(hB_\lambda(\theta_{nh}\omega)+I)X_n,X_n\rangle
+h\langle F(X_n),X_n\rangle\\
&=\langle(hB_\lambda(\theta_{nh}\omega)+I)X_n,X_n\rangle\\
&=(1-h\sigma_\lambda(\theta_{nh}\omega))x_n^2+h(\sigma_\lambda(\theta_{nh}\omega)
+\rho_\lambda(\theta_{nh}\omega))x_ny_n\\
&\quad +(1-h)y_n^2+(1-h\beta_\lambda(\theta_{nh}\omega))z_n^2,
\end{align*}
we have
\begin{align*}
|\langle X_{n+1},X_n\rangle|
&\le(1-h\sigma_\lambda(\theta_n\omega))x_n^2+(1-h)y_n^2+(1-h\beta_\lambda(\theta_n\omega))z_n^2\\
&\quad
+\frac{h}2(\sigma_\lambda(\theta_n\omega)+\rho_\lambda(\theta_n\omega))(x_n^2+y_n^2)\\
&\le
(1-\frac{h(\sigma_\lambda(\theta_n\omega)-\rho_\lambda(\theta_n\omega))}2)x_n^2
+(1-h\beta_\lambda(\theta_n\omega))z_n^2
\\
&\quad +[1-\frac{h}2(2-\sigma_\lambda(\theta_n\omega)-\rho_\lambda(\theta_n\omega))]y_n^2\\
&<\|X_n\|^2
\end{align*}
for small $h$ by the assumption
$0<\rho_\lambda(\omega)<\sigma_\lambda(\omega)\le 1$, where the last
``$<$" holds if $X_n\neq (0,0,0)$. That is, we have shown that
$\|X_{n+1}\|\le \|X_n\|$, which implies that any ball $B_r(0)$
centered at the origin with radius $r>0$ is an isolating
neighborhood for $\varphi_\lambda$, $\lambda\in[0,1]$. It is easy to
verify that $\{(0,0,0)\}$ is an isolated invariant set for the
unperturbed discrete Lorenz system (\ref{unp}) with nontrivial
Conely index. Hence by the continuity property of random conley
index we have that
$h(S_\lambda,\varphi_\lambda)=h(S_0,\varphi_0)\neq\underline 0$,
$\lambda\in[0,1]$. Therefore, by the Waz\.ewski property of random
Conley index we obtain that there exists an invariant random compact
set in $B_r(0)$ for the perturbed discrete Lorenz system. Since the
above argument holds for arbitrary $r>0$, we obtain that
$\{(0,0,0)\}$ is an invariant random compact set for the perturbed discrete Lorenz system.\\
(2) Assume that for each $\lambda\in[0,1]$,
$\sigma_\lambda(\omega),\beta_\lambda(\omega)>0$ almost surely. We
consider the Lorenz system on the sub-manifold:
$M=\{(x,y,z)\in\mathbb R^3|y=0\}$ (clearly $M$ is an invariant
random set of RDS $\varphi_\lambda$, for each $\lambda\in[0,1]$),
then by the computation in (1) we have
\[
\langle X_{n+1},X_n\rangle
=(1-h\sigma_\lambda(\theta_n\omega))x_n^2+(1-h\beta_\lambda(\theta_n\omega))z_n^2.
\]
Hence for small $h$, we have
\[
|\langle X_{n+1},X_n\rangle|\le \|X_n\|^2,
\]
i.e.
\[
\|X_{n+1}\|<\|X_n\|
\]
almost surely whenever $X_n\neq(0,0,0)$. This indicates that when
restricted on the sub-manifold $M$, $\varphi_\lambda$,
$\lambda\in[0,1]$ admits any neighborhood of $(0,0,0)$ in $M$ as an
isolating neighborhood. Clearly the Conley index for $\varphi_0$,
i.e. the unperturbed discrete Lorenz system (\ref{unp}), is
nontrivial. Similar to the argument of (1) we immediately obtain
that $\{(0,0,0)\}$ is an invariant random compact set for the
perturbed discrete Lorenz system when it is restricted on $M$.
\end{example}

\begin{example}\rm
Assume that $\varphi$ is a locally tempered analytic random
diffeomorphism in $\mathbb C^d$ with a fixed point $x=0$ and its
linearization $A(\omega)$ satisfies the conditions of the
multiplicative ergodic theorem. Moreover, if all Lyapunov exponents
$\lambda_i(\omega)$, $1\le i\le p(\omega)$, have the same sign and
are non-resonant, then by (iii) of Main Theorem of \cite{LL} we have
$\varphi$ is analytically conjugate to its linear part. That is, if
we write $\varphi(\omega,x)$ as
\[
\varphi(\omega,x)=A(\omega)x+f(\omega,x),
\]
where $A(\omega)=D\varphi(\omega,0)\in Gl(d,\mathbb C),
f(\omega,0)=0$ and $Df(\omega,0)=0$, then there exists an analytic
random diffeomorphism $h(\omega,x)$ defined in a tempered ball
$V(\omega)$ with $h(\omega,0)=0$ such that
\begin{equation}\label{ex1}
h(\theta\omega,\varphi(\omega,x))=A(\omega)h(\omega,x),~{\rm
for}~x\in V(\omega) ~~{\rm a.e.}
\end{equation}
For the details, the reader can refer to \cite{LL}. Since all
Lyapunov exponents $\lambda_i(\omega)$ have the same sign, it is
clear that $\{0\}$ is an isolated invariant set for $A(\omega)x$ and
$h(\{0\},A(\omega))\neq\underline0$. By (\ref{ex1}) we obtain that
$\{0\}$ is an isolated invariant set for $\varphi$. Further more,
the random Conley index of $\varphi$ is the same as that of linear
part. To see this, consider the family of random diffeomorphisms:
\[
\varphi_\lambda(\omega,x)=A(\omega)x+\lambda
f(\omega,x),~~\lambda\in[0,1].
\]
By (iii) of Main Theorem of \cite{LL} we have for each
$\lambda\in[0,1]$, there exists $V_\lambda(\omega)$ such that
$\varphi_\lambda$ is analytically conjugate to $A(\omega)x$ in
$V_\lambda(\omega)$ and $V_\lambda(\omega)$ play the role of
isolating neighborhood with respect to $\varphi_\lambda$, isolating
$\{0\}$. By a compactness argument there exists a $V(\omega)$ such
that $\varphi_\lambda$ is analytically conjugate to $A(\omega)x$ in
$V(\omega)$ for each $\lambda\in[0,1]$ and $V(\omega)$ is an
isolating neighborhood for each $\varphi_\lambda$. Hence by the
continuity property of random Conely index we have
$h(S_\lambda,\varphi_\lambda)=h(\{0\},\varphi_0)$, where
$S_\lambda(\omega):={\rm Inv}(V(\omega),\varphi_\lambda)$,
$\lambda\in[0,1]$. Since $V(\omega)$ can be chosen sufficiently
small such that $S_1(\omega)=\{0\}$, we obtain that
$h(\{0\},\varphi_1)=h(\{0\},\varphi_0)$, i.e.
$h(\{0\},\varphi)=h(\{0\},A(\omega))$.

If $\varphi$ is a $C^\infty$ locally tempered random diffeomorphism
in $\mathbb R^d$ with a hyperbolic fixed point $x=0$ and its
linearization $A(\omega)$ satisfies the conditions of the
multiplicative ergodic theorem. And all Lyapunov exponents satisfy
non-resonance condition, then by (ii) of Theorem 1.1 of \cite{LL2}
we have $\varphi$ is conjugate to its linear part. By the similar
argument to that of analytic case in $\mathbb C^d$ we can obtain
that $\{0\}$ is an isolated invariant set for $\varphi$ and it has
the same random Conley index with respect to $\varphi$ as it has
with respect to its linear part.
\end{example}

\section*{Acknowledgements}
I am most indebted to my advisor, Professor Yong Li, not only for
his direct helpful suggestions but primarily for his continual
instruction, encouragement and support over all these years. I am
very grateful to Professor Yingfei Yi for helpful discussions and
valuable suggestions during the period when he visit Jilin
University, to Professor Ludwig Arnold for his kind encouragement
and to Professor Dave Richeson for his patient explaining his joint
paper with Professor John Franks.

\end{document}